\documentclass{article}

\usepackage{pgf}
\usepackage{amssymb}
\usepackage{dsfont}
\usepackage{graphicx}
\usepackage{amsmath}
\usepackage{frenchle}
\usepackage[latin1]{inputenc}
\usepackage[T1]{fontenc}

\newtheorem{theo}{Theorem}[section]
\newtheorem{definition}[theo]{Definition}
\newenvironment{defi}{\begin{definition}\rm}{\end{definition}}
\newtheorem{remarque}[theo]{Remark}
\newenvironment{remark}{\begin{remarque}\rm}{\end{remarque}}
\newtheorem{exemple}[theo]{Example}
\newenvironment{ex}{\begin{exemple}\rm}{\end{exemple}}
\newtheorem{lemma}[theo]{Lemma}
\newtheorem{propo}[theo]{Proposition}
\newtheorem{coro}[theo]{Corollary}
\newtheorem{prf}{\it {Proof}}

\newenvironment{demo}{\begin{prf}\rm}{\hfill$\Box$\end{prf}}
\newtheorem{nota}[theo]{Notation}
\newenvironment{notation}{\begin{nota}\rm}{\end{nota}}

\begin{document}

\title{A differential Puiseux theorem in generalised series fields of finite rank.} 

\author{Micka\"{e}l Matusinski}



\maketitle


{\small\noindent \textsc{Address:} \ Universit\"{a}t Konstanz\\
\indent \indent \indent \indent  Fachbereich Mathematik und Statistik\\
\indent \indent \indent \indent 78457 Konstanz, Allemagne.\\
\textsc{email:}\indent \ \ \  mickael.matusin@gmail.com\\
\textsc{Page web:}\ http://sites.google.com/site/mickaelmatusinski/}

\begin{abstract}
 Nous étudions des équations différentielles $F(y,\ldots,y^{(n)})=0$ où $F$ est une séries formelle en $y,y',\ldots,y^{(n)}$, à coefficients dans un \emph{corps de séries généralisées} $\mathds{K}_r$ de rang fini $r\in\mathbb{N}^*$. Notre objet est d'exprimer le support -c'est-à-dire l'ensemble des exposants $\textrm{Supp}\ y_0$- des éléments $y_0\in\mathds{K}_r$ solutions, en fonction des supports des coefficients de l'équation, dont l'union est notée $\textrm{Supp}\ F$.

\begin{center}
\textbf{Abstract}
\end{center} We study differential equations  $F(y,\ldots,y^{(n)})=0$ where $F$ is a formal series in $y,y',\ldots,y^{(n)}$ with coefficients in some field of \emph{generalised power series} $\mathds{K}_r$ with finite rank $r\in\mathbb{N}^*$. Our purpose is to express the support $\textrm{Supp}\ y_0$, i.e. the set of exponents, of the elements $y_0\in\mathds{K}_r$ that are solutions, in terms of the supports of the coefficients of the equation, namely $\textrm{Supp}\ F$.\\

\end{abstract}



\newpage

\tableofcontents

\section{Introduction}
\subsection{About differential Puiseux theorems.} In his PhD thesis \cite[Theorem 12.2]{vdh:autom-asymp}, van der Hoeven proves that:\\
\textsl{Any well-ordered transseries solution to an algebraic differential equation with grid-based transseries coefficients is itself grid-based.}

 Another version of this result can be found in \cite[Corollary 8.38]{vdh:transs_diff_alg}. Without entering into the details of definitions, transseries are formal series built from $\mathbb{R}$, $x$, field operations, composition with exponential and logarithmic functions, and an infinite summation process analogous to the construction of generalised series. Grid-based transseries are those whose support is included in a grid, that is to say the translation of a lattice. They were introduced by J. Ecalle in his proof of Dulac conjecture \cite{ecalle:dulac}. Well-ordered transseries are those built with well-ordered supports. 

This result for transseries addresses two difficulties concerning formal resolution of ordinary differential equations. On the one hand, it shows that what we may call the rank of any well-ordered transseries solution $f$ (i.e. the Archimedean rank of the elements of the support of $f$) is finite, following the fact that it is finite for the grid-based coefficients of the equation \cite[Theorem 4.15: the rank of any grid-based transseries is bounded by the number of elements of a transbasis for it]{vdh:transs_diff_alg}.

 On the other hand, the support of such a finite rank well-ordered transseries solution is actually grid-based. This is a Puiseux type result \cite[for instance]{kr_pa:primer-real-ana-funct} in the context of differential equations, which generalizes several recent results. Grigoriev and  Singer \cite[Corollary 3.1]{grigo:sing} considered polynomial differential equations $P(y,\ldots,y^{(n)})=0$ with $P\in\mathbb{Q}[[x]][y,\ldots,y^{(n)}]$ (i.e. the coefficients are formal power series in $x$ with rational coefficients), and solutions that are power series with \textsl{real} exponents. Then they showed that these exponents belong to a finitely generated $\mathbb{Z}$-module in $\mathbb{R}$. 

Independently, J.Cano \cite[Theorem 1]{jcano:series} studies equations $f(x,y,y'\ldots,y^{(n)})=0$, where
$f\in\mathbb{R}[[x]][[ y,\ldots,y^{(n)}]]$, i.e. $f$ is a power series in the $n+1$ variables $y,\ldots,y^{(n)}$ and in $x$. He shows that any series with \textsl{rational} exponents greater than $n$ which is a solution, is in fact a Puiseux series.

 In their proof of a desingularisation theorem, F.Cano, Moussu and Rolin used a generalisation \cite[Appendix]{cano_moussu_rolin} of J. Cano's result to the case of real exponents:\\
 \textsl{Given a differential equation $f(x,y,xy',xy'')=0$, where
$f\in\mathbb{R}[[x]][[ y,xy',xy'']]$, then the support of any solution is included in a \textsl{lattice} (i.e. a finitely generated sub-semi-group of $\mathbb{R}_{\geq 0}$)}.\\

 The main purpose of this paper is to prove a generalisation of this second Puiseux type part of van der Hoeven's result (and therefore a generalisation of the other cited results) to the context of \textsl{generalised series fields} (which have well-ordered supports and do not carry any log-exp structure a priori) and for \textsl{differential equations defined by formal series} (not only by polynomials: see Equation \ref{eq:equa-diff}). Given a totally ordered abelian group $\Gamma$ and a field $\mathcal{C}$, a generalised series with exponents in $\Gamma$ and coefficients in $\mathcal{C}$ is a formal sum $a=\sum_{\gamma\in \Gamma} a_\gamma t^\gamma$, where $t$ is an abstract variable, the coefficients $a_\gamma=a(\gamma)$ belong to $\mathcal{C}$ and its support $\textrm{Supp}\ a=\{\gamma\in \Gamma \ |\  a_\gamma \neq 0 \}$ is well-ordered in $\Gamma$ (for classical definitions and properties in well-ordering theory, see \cite{bour:theo_ens} and \cite{krivine:theo_ensembles}). The series with empty support is denoted 0. The set of such series endowed with component-wise sum and convolution product is a field \cite{hahn:nichtarchim}, which we denote by $\mathds{K}$. We endow it also with the classical valuation, i.e. the surjective map $v\ :\ \mathds{K}\rightarrow \Gamma\cup\{\infty\}$ given by $v(a)=\min(\textrm{Supp}\ a)$ for $a\neq 0$ and $v(0)=\infty$, with the usual conventions.
 
In particular, we consider a group of exponents $\Gamma$ of \textsl{finite rank} $r\in\mathbb{N}^*$ (i.e. with a finite number of Archimedean classes) and $\mathbb{R}$ as field of coefficients (see Section \ref{sec:defi}). We denote by $\mathds{K}_r$ the corresponding series field. We also suppose that this field is endowed with a derivation which behaves like the one in Hardy fields (see Definition \ref{hypo_deriv}). We consider differential equations as follows:
\begin{defi}\label{defi:equa_diff}
Let 
\begin{equation}\label{eq:equa-diff}
F(y,\ldots,y^{(n)})=0
\end{equation}
be given, where $F$ is a non trivial differential series:
\begin{equation}\label{eq:serie-diff}
F(y,\ldots,y^{(n)})=\sum_{I\in\mathbb{N}^{n+1}}
c_Iy^{(I)}\in\mathds{K}_r[[y,\ldots,y^{(n)}]]\setminus\{0\}.
\end{equation} of arbitrary fixed order $n\in\mathbb{N}$, with coefficients in $\mathds{K}_r$ and \textsl{well-ordered support} $\textrm{Supp}\ F=\cup_{I\in\mathbb{N}^{n+1}}\textrm{Supp}\ c_I$. For any $I\in\mathbb{N}^{n+1}$, we denote $y^{(I)}=y^{i_0}(y')^{i_1}\cdots(y^{(n)})^{i_n}$.
\end{defi}
We suppose that the support of the equation is well-ordered so as to remain in $\mathds{K}_r$ when evaluating the Differential Series \ref{eq:serie-diff} at some series $y_0\in\mathds{K}_r$. Note that this requirement is automatically fulfilled when \ref{eq:serie-diff} is a polynomial.

\begin{ex}.\\
1. There are natural settings in which such differential fields of series arise: algebraic, analytic or even formal vector fields. For example, following \cite{cano_moussu_rolin}, take $X\ :\ t_k'=F_k(t_1,\ldots,t_r),\ k\in\{1,\ldots,r\}$ to be an analytic vector field over a real analytic $r$-dimensional manifold $m$. Consider $\gamma : t\mapsto \gamma(t)$, $t\geq 0$, an integral curve of X having a unique $\omega$-limit point $p$ and assume that $\gamma$ is sub-analytically non-oscillating. It is shown in Section 2 of \cite{cano_moussu_rolin} that one can associate a Hardy field $K_\gamma$ to $\gamma$ which has at most rank $r$. So, in the case of maximal rank $r$, we can take the formal counterpart of such a Hardy field, namely the field $\mathds{K}_r$ endowed with the derivation determined by $X$.\\
more generally, now take $X\ :\ t_k'=F_k(t_1,\ldots,t_r),\ k\in\{1,\ldots,r\}$ to be a formal vector field. For all $i\in\{1,\ldots,r\}$, denote $\tau^{(k)}=(\tau_1^{(k)},\ldots,\tau_r^{(k)})$ to be the least multi-exponent for the lexicographical ordering among those appearing in the monomials $t_1^{\tau_1^{(k)}}\cdots t_r^{\tau_r^{(k)}}$ of $F_k(t_1,\ldots,t_r)$. Lemma \ref{prop:val_d_k} gives the following:\\
\textbf{Claim}.\textsl{ The corresponding field of series $\mathds{K}_r$ can be endowed with a Hardy type derivation if and only if the matrix $(\tau_j^{(k)})_{1\leq j,k\leq r}$ is such that:
\begin{itemize}
    \item for any $j=1,\cdots,k-1$, $\ \tau_j^{(k+1)}=\tau_j^{(k)}$;
\item $\tau_k^{(k+1)}=\tau_k^{(k)}-1$;
\item $(0,\ldots,0,\tau_{k+1}^{(k+1)}-1,\ldots,\tau_r^{(k+1)})>_{lex} (0,\ldots,0,\tau_{k+1}^{(k)},\ldots,\tau_r^{(k)}).$
\end{itemize}}
It may be interesting to classify, up to a change of coordinates, the vector fields which verify such a property.\\
\noindent 2. As cited before, one can take any finite rank differential subfield of the field of well-ordered transseries. For instance, take the field of generalised series $\mathbb{R}((\Gamma))$ with $\Gamma=\{\log(x),x,\exp(x)\}^\mathbb{R}$, the group of words in  $t_1=\exp(-x)$, $t_2=1/x$ and $t_3=1/\log(x)$ with real exponents, the usual comparison relations and derivation. Thus $t_1'=-t_1$, $t_2'=-t_2^2$ and $t_3'=-t_1t_2^2$.
\end{ex}  

\subsection{Content of our article.} One of the main results of this paper is:
\begin{theo}\label{theo_princ1}
Given a Differential Equation \ref{eq:equa-diff}, the support of any solution $y_0\in\mathds{K}_r$ such that $v(y_0^{(i)})>\underline{0}$ for any $i=0,\ldots,n$, is obtained by finitely many elementary transformations from the supports of $F$ and the $t_k'/t_k$'s, $k=1,\ldots,r$. 
\end{theo}
The reader will find the necessary precisions in Definition \ref{hypo_deriv} and Notation \ref{rem_decompo}. We emphasise that the $t_k'/t_k$'s do not depend on the Equation \ref{eq:equa-diff}, but only on the differential field $\mathds{K}$ that we consider. This result is obtained as an immediate corollary of the Theorem \ref{theo_princ}.

Our notion of ``being deduced by finitely many elementary transformations" generalises those of ``being grid-based" and ``belonging to a lattice" to the setting of well-ordered supports (see Remark \ref{lemme:reseaux}). Moreover, beyond this Puiseux type result, our Theorem \ref{theo_princ} shows a dichotomy result identical to the classical one for pseudo-Cauchy sequences in valuation theory \cite[Section 0.1]{kuhl:ord-exp}. In the case of non resolution of the equation, we have a situation of stabilisation of the valuation, which is analogous to the one of monomialization proved in \cite{matu_rolin:subanaED} for sub-analytic differential equations.

The assumption that  $v(y_0^{(i)})>\underline{0}$ for any $i=0,\ldots,n$, means that $y_0$ and all its derivatives are infinitesimal. In other words, we can say that the solution is supposed to ``tend to $(0,\ldots,0)\in\mathbb{R}^{n+1}$", this point being possibly a singular point for the Equation \ref{eq:equa-diff}. We intend to obtain the same result for arbitrary generalised series solutions, i.e. possibly non infinitesimal. In the Section \ref{sect:posit-val}, we show that such is the case for polynomial differential equations: see Theorem \ref{theo:poly-diff}.\\

Our paper is organised into four parts. In Section \ref{sect:main_result}, we state the main Theorem \ref{theo_princ}. Since we are working with formal equations as \ref{eq:equa-diff} rather than with polynomial ones, there may be some restrictions on the operations we can make. In Section \ref{sec:defi}, we check whether the Differential Series \ref{eq:serie-diff} we consider can be evaluated at some series $y_0\in\mathds{K}_r$. In doing so, we also characterise the support of the obtained generalised series $F(y_0,\ldots,y_0^{(n)})$. Section \ref{sec:transfo-series-diff} deals with three transformations of differential series. Two of these transformations, namely \textsl{additive conjugations} and \textsl{multiplicative conjugations}, are the ones used already in \cite{vdh:autom-asymp} and \cite{vdh:transs_diff_alg}. The third and most difficult ones, namely changes of derivations, are not needed in van der Hoeven's work since he uses upward and downward shiftings corresponding to the logarithmic-exponential structure. A first application of these transformations, in Section \ref{sect:posit-val}, is the reduction of the case of polynomial differential equations with arbitrary generalised series solutions, to our main Theorem \ref{theo_princ}. Then we use them to prove this Theorem \ref{theo_princ}  in Section \ref{sect:weier=w}. This proof, as is the case for the cited results, uses an inductive method, evaluating the Differential Series \ref{eq:serie-diff} at some longer and longer initial parts of a solution. But all the proofs of the above cited differential results, rely on an adaptation of the classical Newton polygonal method to the differential case: the Newton-Fine polygonal method \cite{fine:polygon}. Instead of such a polygonal representation, we use a valuative approach and express directly the relation between the exponents of the solution and those of the coefficients of the equation.

The major part of this work was done while the author was a doctoral fellow with J.-P. Rolin at the University of Bourgogne (see \cite{matu:these}). Another version of this article was written while the author was a visiting post doctoral fellow at the University of Saskatchewan, and was partially supported by Salma Kuhlmann's NSERC Discovery Grant. The author thanks gratefully Salma Kuhlmann for providing many good advises and interesting comments on this work.

\section{The main result}\label{sect:main_result}
\subsection{The Hahn field of rank \textsl{r}.} For any positive integer $r$, the \textsl{Hahn group of rank $r$} is the product of $r$ copies of $\mathbb{R}$ ordered lexicographically, say $\overrightarrow{\mathbb{R}}^r$. By Hahn's embedding theorem in \cite{hahn:nichtarchim}, any totally ordered group $\Gamma$ with finite rank (i.e. with a finite number of Archimedean classes) embeds in such a Hahn group. From now on, we fix some positive integer $r$ and, \textsl{without loss of generality, we reduce to the corresponding Hahn group $\Gamma=\overrightarrow{\mathbb{R}}^r$} as a group of exponents for the generalised series. We write $\underline{0}$ its neutral element $(0,\ldots,0)$.

We denote by $\mathds{K}_r$ the corresponding field of generalised series \textsl{with real coefficients} and by $\mathds{K}_r^\prec$ (respectively
$\mathds{K}_r^0$, $\mathds{K}_r^\succ$, $\mathds{K}_r^\preccurlyeq$) the subset of $\mathds{K}_r$ given by $\{a\in\mathds{K}_r \ | \
v(a)>0 \ (respectively \ =0, \ <0,\ \geq 0)\}$. The symbols $\prec,\preceq$ and $\succ$ denote the usual dominance relations for functions ($t^\alpha\prec t^\beta\Leftrightarrow\alpha>\beta$). With this notation,  $\mathds{K}_r^\preccurlyeq$ is the valuation ring of $\mathds{K}$ and  $\mathds{K}_r^\prec$ its maximal ideal. We require the coefficients to be real so that they will be compatible with the real exponents of the monomials, applying the Leibniz rule (HD0) below.

We will also use the notion of \textsl{leading term} of a series $a$, namely  $\delta(a)=a_{v(a)}t^{v(a)}$ for any non zero $a\in\mathds{K}_r$, and the classical \textsl{equivalence relations}\begin{center}
 $\begin{array}{lcl}a\asymp b&\Leftrightarrow& v(a)=v(b)\\
 a\sim b&\Leftrightarrow& v(a-b)>\min\{v(a),v(b)\}\end{array}$.
\end{center}

From now on, we also use another notation for elements of $\mathds{K}_r$. For any $i\in\{1,\ldots,r\}$, we set $t_k=t^{e_k}$, where
$e_k=(0,\ldots,1,\ldots,0)$, with 1 in the $k^{th}$ position is the $k^{th}$ generator of the $\mathbb{R}$-vector space $\mathbb{R}^r$. So any element $a\in\mathds{K}_r$ can be written \begin{center}
$a=\sum_{\alpha\in\Gamma}a_\alpha t^\alpha=\sum_{\alpha_1\in\mathbb{R}}t_1^{\alpha_1} \left(\sum_{\alpha_2\in
I_{\alpha_1}}t_2^{\alpha_2}\left(\cdots \left(\sum_{\alpha_r\in
I_{\alpha_1,\ldots,\alpha_{r-1}}} a_\alpha
t_r^{\alpha_r}\right)\cdots\right)\right)$,
\end{center} where the $I_{\alpha_1,\ldots,\alpha_k}$'s are subsets of $\mathbb{R}$. In particular, any element $a$ of $\mathds{K}_r^\prec$ can be written\begin{center}
$\begin{array}{lcl}
 a&=& \sum_{\alpha_r>0}a_{(0,\ldots,0,\alpha_r)}t_r^{\alpha_r}+
 \sum_{\alpha_{r-1}>0}t_{r-1}^{\alpha_{r-1}} \left(\sum_{\alpha_r\in
I_r}a_{(0,\ldots,0,\alpha_{r-1},\alpha_r)}t_r^{\alpha_r}\right)+\cdots\\
&=& a_{r}+\cdots+a_{1},
\end{array}$
\end{center}
with $v(a_i)>v(a_j)$ for any $i,j\in\{1,\ldots,r\}$, $i<j$, whenever $a_j\neq 0$.

\begin{notation}\label{rem_decompo}
For any $k\in\{1,\ldots,r\}$, we denote by $\mathds{K}_{r,k}$ the additive subgroup of $\mathds{K}_r^\prec$, of elements \begin{center}
$a_k= \sum_{\alpha_{k}>0}t_{k}^{\alpha_{k}}  \left(\sum_{\alpha_{k+1}\in
I_{\alpha_{k}}}t_{k+1}^{\alpha_{k+1}}\left( \cdots  \left(\sum_{\alpha_r\in
I_{\alpha_k,\ldots,\alpha_{r-1}}} a_{(0,\ldots,0,\alpha_{k},\ldots,\alpha_r)}t_r^{\alpha_r} \right)\cdots\right)\right)$, 
\end{center}
together with the series 0. 
Note that  $\mathds{K}_r^\prec=\mathds{K}_{r,r}\oplus\cdots\oplus\mathds{K}_{r,1}$  with $\underline{0}<v(\mathds{K}_{r,r}^*)<\cdots<v(\mathds{K}_{r,1}^*)$.
\end{notation} 

\subsection{Endowing the series with a Hardy type derivation.} We suppose from now on that $\mathds{K}_r$ is endowed with a derivation $D_0$ as follows:
\begin{defi}[Hardy type derivation]\label{hypo_deriv} We say that a map $D_0:\mathds{K}_r\rightarrow \mathds{K}_r$, $a\mapsto a'$ with $t_k'\neq 0$ for any $k=1,\ldots,r$, is a \emph{Hardy type derivation} if the following
conditions are satisfied:
\begin{description}
    \item[(HD0)] $\forall \alpha\in\Gamma$, $D_0(t^\alpha)=(t^\alpha)'=(t_1^{\alpha_1}\cdots
t_r^{\alpha_r})'=t^{\alpha}(\alpha_1t'_1/t_1+
\cdots+\alpha_rt'_r/t_r)$;
    \item[(HD1)]  for any $a=\sum_{\alpha\in\Gamma}a_\alpha
t^\alpha$ in $\mathds{K}_r$, $D_0(a)=a'=\sum a_\alpha (t^\alpha)'$;
    \item[(HD2)]   for all $k=1,\ldots,r-1$, $v(t_k')> v(t_{k+1}')$;
    \item[(HD3)]  for all $k=1,\ldots,r-1$, $v(t'_k/t_k)<
v(t'_{k+1}/t_{k+1})$.
\end{description}
For any $k\in\{1,\ldots,r\}$, we denote $\delta(t'_k)=T_kt^{\tau^{(k)}}$, $d_k:=t^{\tau^{(k)}}/t_k$.
\end{defi}

\begin{remark}\label{rem:Hardy_deriv}
 Such a Hardy type derivation $D_0$ is a derivation in the usual sense, i.e. for any $a,b\in\mathds{K}$, $(a+b)'=a'+b'$ and $(ab)'=a'b+ab'$. We note also that such a derivation is uniquely determined by its restriction to the $t_k$'s. Indeed, for any non zero $a\in\mathds{K}_r$:
\begin{center}
$\begin{array}{lcl}
a'&=&=\sum a_\alpha (t^\alpha)'\\
&=&\sum a_\alpha \alpha(\alpha_1t_1'/t_1+\cdots+\alpha_rt_r'/t_r) \\
&=&[\sum (a_\alpha\alpha_1) \alpha] t_1'/t_1\ +\cdots+[\sum(a_\alpha\alpha_r) \alpha] t_r'/t_r.
\end{array}$ 
\end{center}

\end{remark}

Hypothesis (HD0) is an extension to real exponents of the usual Leibniz rule that holds for rational exponents for any derivation. Likewise, Hypothesis (HD1) is an extension of the linearity property for derivations. Thus $D_0$ is a \textsl{strong derivation} in the sense of van der Hoeven \cite[Section 5.1]{vdh:transs_diff_alg}. 

Hypothesis (HD2) and (HD3) imply that the following crucial properties hold:\\
\textbf{Claim 1}. For any $ a,b\in\mathds{K}_r$:
\begin{description}
\item[(HD4)] \textbf{(L'Hospital's rule)}: \textrm{if } $v(a),v(b)\neq 0$, then \begin{center}
$v(a)\leq v(b)\Leftrightarrow v(a')\leq v(b')$;
\end{center} 
\item[(HD5)] \textbf{(Compatibility with the logarithmic derivatives)}: \begin{center}
$|v(a)|\geq |v(b)|>0 \Rightarrow v(a'/a)\leq v(b'/b)$.\\
Moreover, $v(a'/a)= v(b'/b)\Leftrightarrow v(a),v(b)\textrm{ are Archimedean equivalents}$.
\end{center}
\end{description}
\begin{demo}
By (HD1), to prove l'Hospital's rule, it suffices to check it for monomials. Thus, consider $t^\alpha,t^\beta$ with $\alpha\geq\beta$. By (HD0), $(t^\alpha)'=t^\alpha(\alpha_kt_i'/t_i+\cdots+\alpha_rt_r'/t_r)$ and $(t^\beta)'=t^\beta(\beta_jt_j'/t_j+\cdots+\alpha_rt_r'/t_r)$ where $\alpha=(0,\ldots,0,\alpha_i,\cdots,\alpha_r)$ and $\beta=(0,\ldots,0,\beta_j,\ldots,\beta_r)$. Using (HD2) and (HD3), it is a routine to verify that $v((t^\alpha)')=\alpha+d_i\geq v((t^\beta)')=\beta+d_j$.

Consider now $ a,b\in\mathds{K}_r$ such that $|v(a)|=\alpha\geq |v(b)|=\beta>0$. Suppose that $\alpha=(0,\ldots,0,\alpha_i,\cdots,\alpha_r)$ and $\beta=(0,\ldots,0,\beta_j,\ldots,\beta_r)$. By (HD1) and (HD5), $v(a')=\alpha+d_i$ and $v(b')=\beta+d_j$. So, $v(a'/a)=d_i$ and $v(b'/b)=d_j$. The Property (HD5) therefore follows from (HD3).
\end{demo}

\begin{remark}\label{rem:H-fields}
In our context of a Hardy type derivation, the valuation $v$ is a \textsl{differential valuation} in the sense of Rosenlicht \cite{rosenlicht:diff_val}. Moreover, 
(HD5) is an additional property that holds in Hardy fields \cite[Propositions 3 and 4]{rosenlicht:rank} as well as in pre-\textsl{H}-fields \cite[Lemma 3.5]{vdd:aschenbrenner:asymptotic}. Otherwise we are in the case that Rosenlicht calls a \textsl{Hardy type valuation} \cite{rosenlicht:val_gr_diff_val2}. 

To obtain a structure of \textsl{H-field} \cite{vdd:aschenbrenner:asymptotic}, it suffices to endow $\mathds{K}_r$ with the ordering associated to the notion of \textsl{leading coefficient}, and to require that for any $k$, $t_k'>0$, i.e. $T_k>0$.
\end{remark}
We note also that:\\
\textbf{Claim 2}. The field of constants of $D_0$ is $\mathbb{R}\subset\mathds{K}_r$. 
\begin{demo}If there was a non constant series $a$ with $a'=0$, then by (HD4), this would imply that there exists a non constant monomial $m$ with derivative 0. Then, taking any other monomial $\tilde{m}$ with non zero derivative, we would have either $v(\tilde{m})<v(m)$ or $-v(\tilde{m})<v(m)$, but $v(m')=\infty<v((\tilde{m}^\pm)')$. This contradicts (HD4).
\end{demo} 

With S. Kuhlmann in \cite{matu-kuhlm:hardy-deriv-EL-series}, we develop these ideas in greater details for more general series fields. To conclude this section, we give some key properties of the values $v(t_k')=\tau^{(k)}$. 
\begin{nota}\label{nota:theta}
  For any $k=1,\ldots,r$, if $v(d_k)\neq \underline{0}$, then we write, \begin{center}
$v(d_k)=\theta^{(k)}= (0,\ldots,0,\theta^{(k)}_{\tilde{k}},\ldots,\theta^{(k)}_r)$.
\end{center}
with $\theta^{(k)}_{\tilde{k}}\neq 0$.
\end{nota}

\begin{propo}\label{prop:val_d_k}
For any $k>l$, there exists $m\in\{l,\ldots,r\}$ such that: \begin{center}
$v(d_k/d_l)=(0,\ldots,0,\theta_{m}^{(k)}- \theta_{m}^{(l)},\ldots,\theta_{r}^{(k)}- \theta_{r}^{(l)})$ with $\theta_{m}^{(k)}- \theta_{m}^{(l)}>0$.
\end{center}
\end{propo}
\begin{demo} We show that:\\
$\indent\ \bullet\ $ for any $i=1\cdots k-1$, $\tau_i^{(k+1)}=\tau_i^{(k)}\ ; \\
\indent\ \bullet\ \tau_k^{(k+1)}=\tau_k^{(k)}-1\ ;\\
\indent\ \bullet\ (0,\ldots,0,\tau_{k+1}^{(k+1)}-1,\ldots,\tau_r^{(k+1)})> (0,\ldots,0,\tau_{k+1}^{(k)},\ldots,\tau_r^{(k)}).$\\
From Definition \ref{hypo_deriv}, $v(t'_k)>v(t'_{k+1})$ and $v(d_k)<v(d_{k+1})$. So \begin{center}
$(\tau_{1}^{(k+1)}-\tau_{1}^{(k)},\ldots,\tau_{k}^{(k+1)}-\tau_{k}^{(k)}, \tau_{k+1}^{(k+1)}-\tau_{k+1}^{(k)},\ldots,\tau_{r}^{(k+1)}-\tau_{r}^{(k)}) <\underline{0}$
\end{center} and \begin{center}
$(\tau_{1}^{(k+1)}-\tau_{1}^{(k)},\ldots,\tau_{k}^{(k+1)}-\tau_{k}^{(k)}+1, \tau_{k+1}^{(k+1)}-1-\tau_{k+1}^{(k)},\ldots,\tau_{r}^{(k+1)}-\tau_{r}^{(k)}) >\underline{0}$.
\end{center} Thus, for any $i=1\cdots k-1$, $\tau_i^{(k+1)}=\tau_i^{(k)}$ and $\tau_k^{(k)}\geq\tau_k^{(k+1)}\geq\tau_k^{(k)}-1$. But, for all $\alpha_k,\beta_{k+1}\in\mathbb{R}^*$, $(t_k^{\alpha_k})'=\alpha_kt_k^{\alpha_k}t'_k/t_k$ and $(t_{k+1}^{\beta_{k+1}})'=\beta_{k+1}t_{k+1}^{\beta_{k+1}} t'_{k+1}/t_{k+1}$. Moreover, for any $\alpha_k>0$, $v(t_k^{\alpha_k})>v(t_{k+1}^{\beta_{k+1}})$ which implies that $v((t_k^{\alpha_k})')>v((t_{k+1}^{\beta_{k+1}})')$. So for any $\alpha_k>0$ and $\beta_{k+1}\in\mathbb{R}^*$, we obtain that \begin{center}
$(0,\ldots,0,\alpha_k, -\beta_{k+1}, 0,\ldots,0)> v(d_{k+1}/d_k) = (0,\ldots,0,\tau_k^{(k+1)}-\tau_k^{(k)}+1,\tau_{k+1}^{(k+1)}-1-\tau_{k+1}^{(k)}, \ldots, \tau_{r}^{(k+1)}-\tau_{r}^{(k)})$,
\end{center} 
and so \begin{center}
$\tau_k^{(k+1)}-\tau_k^{(k)}+1\leq 0\Leftrightarrow \tau_k^{(k+1)}\leq\tau_k^{(k)}-1$.
\end{center} Thus $\tau_k^{(k+1)}=\tau_k^{(k)}-1$. 
We conclude recalling that $v(d_k)<v(d_{k+1})$.
\end{demo}

\begin{coro}\label{coro:k_0}
With the Notation \ref{nota:theta}, the following dichotomy holds:\begin{itemize}
    \item either there exists $k\in\{1,\ldots,r\}$ such that $\tilde{k}=k$. Then such $k$ is unique. We denote it by $k_0$.
Moreover, for any $k<k_0$, $\tilde{k}>k$, and for any $k\geq k_0$, $\theta^{(k)}$ and $\theta^{(k_0)}$;
\item or, for any $k=1,\ldots,r-1$, $\tilde{k}>k$, and $\theta^{(r)}=\underline{0}$. Then we write $k_0=r$.
\end{itemize}
\end{coro}

\begin{coro}\label{coro:deriv-val-0}
$v((t^\alpha)')=\underline{0}$ if and only if $\alpha=-v(d_{k_0})$ or $\alpha=\underline{0}$.
\end{coro}

\begin{coro}\label{coro:deriv_neg}
Consider $k_0$ as defined in the Corollary \ref{coro:k_0}, and $k\in\{1,\ldots,r\}$. If $v(d_k)=\underline{0}$, then $d_k^{(i)}=0$ for any $i\in\mathbb{N}^*$. If not, then we have:
\begin{enumerate}
\item $\textrm{ if }k\geq k_0\textrm{, then } v(d_k^{(i)})=v(d_k)+i v(d_{k_0})$ for any $i\in\mathbb{N}$;
\item $\textrm{ if }k< k_0\textrm{, then}$:\begin{enumerate}
\item if $\tilde{k}<k_0$, then $v(d_k^{(i)})=v(d_k)+iv(d_{\tilde{k}})$  for any $i\in\mathbb{N}^*$;
\item if $\tilde{k}>k_0$, then $v(d_k^{(i)})=v(d_k)+v(d_{\tilde{k}})+(i-1)v(d_{k_0})$ for any $i\in\mathbb{N}^*$;
\item if $\tilde{k}=k_0$, then:\begin{itemize}
\item either $\exists j\in\mathbb{N}^*$, $\theta^{(k)}_{k_0}=-j\theta^{(k_0)}_{k_0}$. Then,  $v(d_k^{(i)})=v(d_k)+iv(d_{k_0})$ for any $i=1,\ldots,j$, and $v(d_k^{(i)})=v(d_k)+v(d_{\hat{k}})+(i-1)v(d_{k_0})$ for some $\hat{k}>k_0$ and any $i>j$;
\item or $v(d_k^{(i)})=v(d_k)+iv(d_{k_0})$ for any $i\in\mathbb{N}^*$. 
\end{itemize}
\end{enumerate}
\end{enumerate}
\end{coro}

\subsection{The main Theorem.} To formulate our main theorem, we need the following definitions:
\begin{defi}\label{defi:stabil}
Given a Differential Series \ref{eq:serie-diff} and a non zero generalised series $y_0\in\mathds{K}_r$, we say that:
\begin{itemize}
    \item the series $y_0$ is \textsl{compatible} with \ref{eq:serie-diff} if the family $(c_Iy_0^{(I)})_{I\in\mathbb{N}^{n+1}}$ is strongly summable (see Definition \ref{defi:sommable}). This implies that for any initial part $p$ of $y_0$, the evaluation of $F$ at $p$, is well-defined: $F(p,\ldots,p^{(n)})\in \mathds{K}_r$. A \textsl{solution} of the corresponding Differential Equation \ref{eq:equa-diff} is a compatible $y_0$ such that $F(y_0,\ldots,y_0^{(n)})=0\in \mathds{K}_r$;
\item the series $y_0$, supposed to be compatible with \ref{eq:serie-diff}, \textsl{stabilises} on \ref{eq:serie-diff} with initial part $p_0$ if there exists a proper initial part $p_0$ of $y_0$ such that, for any longer initial part $p$ of $y_0$ (in particular for $p=y_0$), we have $v(F(p,\ldots,p^{(n)}))=v(F(p_0,\ldots,p_0^{(n)}))\neq v(0)$. 
\end{itemize}\end{defi}
Note that, in the case where \ref{eq:serie-diff} is a \textsl{differential polynomial}, any $y_0\in\mathds{K}_r$ is compatible with it. 

 This notion of stabilisation is identical to that of monomialisation of sub-analytic differential equations proved in \cite{matu_rolin:subanaED}.
\begin{defi}\label{elem_transfo} Given two well-ordered subsets $X_1$ and $X_2$ of $\Gamma_{\geq 0}$, elements $\alpha>0$ and $\beta$ of $\Gamma$, we call \textsl{elementary transformations}:
\begin{itemize}
    \item the sum of two sets: $X_1+X_2=\{\xi_1+\xi_2\ |\ \xi_1\in X_1,\xi_2\in X_2\}$;
    \item the generation of the additive semi-group:\\ $\indent\ \left\langle X_1 \right\rangle=\{k_1\xi_1+\cdots+k_q\xi_q\ |\ k_i\in\mathbb{N},\ \xi_i\in X_1,\ q\in\mathbb{N}\}$;
    \item the addition of a new generator: $X_1+\mathbb{N}\alpha$;
    \item the negative translation by $\beta$: $(X_1)_{\geq\beta}-\beta =\{\alpha\in X_1\ |\ \alpha\geq\beta\}-\beta$.
\end{itemize}
\end{defi}

\begin{theo}\label{theo_princ}
Given a Differential Series \ref{eq:serie-diff} and a series $y_0\in\mathds{K}_r$ with  $v(y_0^{(i)})>\underline{0}$ for any $i=0,\ldots,n$ (so, in particular, $y_0$ is compatible with \ref{eq:serie-diff}: see Proposition \ref{compo_cond1}), there exists a well-ordered subset $\mathcal{R}$ of $\Gamma_{>0}$ obtained from $\textrm{Supp}\ F$ and $\textrm{Supp}\ t_k'/t_k$, $k=1,\ldots,r$, by a finite number of elementary transformations such that:
\begin{itemize}
    \item either $\textrm{Supp}\ y_0\subseteq \mathcal{R}$;
\item or the series $y_0$ stabilises on $F$ with initial part $p_0$ and $\textrm{Supp}\ p_0\subseteq \mathcal{R}$. In this case, $y_0$ can not be a solution of the corresponding Differential Equation \ref{eq:equa-diff}.
\end{itemize}
\end{theo}

As an immediate corollary, we obtain Theorem \ref{theo_princ1}. 


It can happen that $y_0$ does not stabilise on $F$ while also $F (y_0)\neq 0$ (for instance, consider $x$ near 0 in $\mathbb{R}_{>0}$, the equation $y-\sum_{k\in\mathbb{N}}x^{k}+\exp(-1/x)=0$ and the series $y_0=\sum_{k\in\mathbb{N}}x^{k}$ with the usual valuation).

\begin{remark}\label{lemme:reseaux} 
We can deduce from the theorem more particular results in the case that the supports are assumed to be grid-based or included in a lattice of $\Gamma_{\geq\underline{0}}$. For instance, suppose that the supports of the equation and of $t_k'/t_k$, $k=1,\ldots,r$, are included in some lattice. It suffices to show that applying each elementary transformation in Definition \ref{elem_transfo} to some lattice again produces a lattice. For the three first transformations, it is a consequence of Proposition 2.1 and Exercise 2.1 in \cite{vdh:transs_diff_alg}. For the fourth one, namely the negative translation, we prove the following claim:\\

\noindent \textbf{Claim}. \emph{Let $l\in\mathbb{N}^*$, $\lambda_k\in\Gamma_{\geq 0}$ for $k\in\{1,\ldots,l\}$, and $\Lambda=\left\langle
\lambda_{1},\ldots,\lambda_{r}\right\rangle $ be the corresponding lattice. Then, for any $\beta\in\Gamma$, $\Lambda_{\geq\beta}-\beta$, the negative translation by $\beta$ of $\Lambda$, is included in some lattice $\left\langle
\nu_{1},\ldots,\nu_{m}\right\rangle $}.
\begin{demo}
We proceed by induction on $l$, the number of generators of $\Lambda$. If $l=1$, we denote $k_{0}:=\min\{ k\in\mathbb{N\,}|\, k\lambda_{1}\geq\beta\} $. Any element $\alpha$ of $\Lambda_{\geq\beta}-\beta=\left\langle \lambda_{1}\right\rangle _{\ge\beta}-\beta$ can be written $k\lambda_{1}-\beta$ with $k\in\mathbb{N}$ and
$k\lambda_{1}\geq\beta$. Thus $\alpha=(k-k_{0})\lambda_{1}+(k_{0}\lambda_{1}-\alpha)$. Then we set $\nu_{1}=\lambda_{1}$ and
$\nu_{2}=k_{0}\lambda_{1}-\beta$.

 If $l\geq 2$, we suppose that the lemma holds for lattices generated by at most $l-1$ elements. We consider a lattice $\Lambda=\left\langle\lambda_{1},\ldots\lambda_{l}\right\rangle $. Let $k_{0}$ be the least natural number $k$ such that $k\lambda_{l}\geq\beta$. Any element $\alpha$ of $\left\langle \Lambda\right\rangle _{\geq\beta}-\beta$ can be written as $\sum_{i=1}^{l-1} k_{i}\lambda_{i}+k_{l}\lambda_{l}-\beta $. If $k_{l}\geq k_{0}$, then \begin{center}
$\alpha=\sum_{i=1}^{l-1}k_{i}\lambda_{i}+(k_{l}-k_{0}) \lambda_{l}+k_{0}\lambda_{l}-\beta$.
\end{center} So $\alpha$ belongs to $\left\langle
\lambda_{1},\ldots,\lambda_{l-1},\lambda_{l},k_{0}\lambda_{l}-\beta\right\rangle$. If $k_{l}< k_{0}$, then $\alpha$ belongs to a set $\left\langle
\lambda_{1},\ldots,\lambda_{l-1}\right\rangle
_{\geq(\beta-k_{l}\lambda_{l})}-(\beta-k_{l}\lambda_{l})$
 to which we apply the induction hypothesis. It follows that $\alpha$ belongs to some lattice $\tilde{\Lambda}_{k_l}$. Then \begin{center}
$\left\langle \Lambda\right\rangle _{\geq\beta}-\beta\subset(\left\langle
\lambda_{1},\ldots,\lambda_{l-1},\lambda_{l},k_{0}\lambda_{l}-\beta\right\rangle+ \sum_{k_l<k_0}\tilde{\Lambda}_{k_l})$.
\end{center}
\end{demo}
\end{remark}

\section{Differential series.}\label{sec:defi}
\subsection{Introducing new derivations.}
From $D_0$ we build $r$ other derivations corresponding to the $r$ Archimedean classes of the value group of $\mathds{K}_r$. 
\begin{defi}\label{defi_D_k}
For any $k\in\{1,\ldots,r\}$, we set $D_k(y)=y'/d_k$
where $d_k=\delta(t'_k/t_k)$ is the leading term of $t'_k/t_k$. 
\end{defi} 
\begin{propo}\label{propo_D_k} For any $k\in\{1,\ldots,r\}$, $D_k$ is a Hardy type derivation on $\mathds{K}_r$ such that, for any monomial $t^\mu=t_k^{\mu_k}\cdots t_r^{\mu_r}$ in $\mathds{K}_{r,k}$, we have \begin{center}
$D_k^i(t_k^{\mu_k}\cdots t_r^{\mu_r}) \sim \mu_k^it_k^{\mu_k}\cdots
t_r^{\mu_r}$.
\end{center}
\end{propo}
\begin{demo} The derivations $D_k$, $k=1,\ldots,r$ inherit the properties of $D_0$ (see Definition \ref{hypo_deriv}). Moreover, by (HD0) for $D_0$, we deduce that for any $\mu=(0,\ldots,0,\mu_k,\ldots,\mu_r)$, \begin{center}
$(t^\mu)'=t^\mu(\mu_k t'_k/t_k+\cdots+\mu_r t'_r/t_r)\sim t^\mu\mu_k t'_k/t_k$,
\end{center}
and so $D_k(t^\mu)\sim \mu_kt^\mu$.
\end{demo}

For our proof of the Theorem \ref{theo_princ}, we will work with differential equations like \ref{eq:equa-diff}, but using as derivation some of these $D_k$, $k=1,\ldots,r$, instead of the initial derivation $D_0$:
\begin{defi}\label{defi:equa_diff-D_k}
Given some $k=1,\ldots,r$, we consider equations
\begin{equation}\label{eq:equa-diff-D_k}
F(y,\ldots,D_k^ny)=0
\end{equation}
where $F$ is a non trivial differential series:
\begin{equation}\label{eq:serie-diff-D_k}
F(y,\ldots,D_k^ny)=\sum_{I\in\mathbb{N}^{n+1}}
c_Iy^{(I)}\in\mathds{K}_r[[y,\ldots,D_k^ny]]\setminus\{0\}.
\end{equation} of arbitrary fixed order $n\in\mathbb{N}$, with coefficients in $\mathds{K}_r$ and \textsl{well-ordered support} $\textrm{Supp}\ F=\cup_{I\in\mathbb{N}^{n+1}}\textrm{Supp}\ c_I$.
\end{defi}
 The problem of checking whether a change of derivation from an Equation \ref{eq:equa-diff} to an Equation \ref{eq:equa-diff-D_k} is well-defined, is the object of the Section \ref{sect:chgement-deriv}.

In the following proposition, we show that for any of these new derivations, the condition ``$v(D_k^i(y_0))>\underline{0}$ for all $i\in\{0,\ldots,n\}$" of the Theorem \ref{theo_princ} reduces to ``$v(y_0)>\underline{0}$".
\begin{propo}\label{deriv_posit}
Let a non zero series $y_0\in\mathds{K}_r$ and $k\in\{1,\ldots,r\}$ be given. If $v(y_0)>\underline{0}$, then $v(D_k^i(y_0))>\underline{0}$  for all $i\in\mathbb{N}$.
\end{propo}
\begin{demo}
 We denote $v(y_0)=\mu=(0,\ldots,0,\mu_l,\ldots,\mu_r)$ for some $l\in\{1,\ldots,r\}$ with $\mu_l\neq 0$. For the case $l=k$, the result follows from Proposition \ref{propo_D_k}.

For the case $k<l$, we show by induction on $i$ that for any $i\in\{0,\ldots,n\}$, \begin{center}
$v(D_k^iy_0)=(0,\ldots,0,\beta_{k_i},\ldots,\beta_r)$
\end{center} for some $k_i>k$ and $\beta_{k_i}>0$. If $i=0$,  \begin{center}
$v(y_0)=(0,\ldots,0,\mu_l,\ldots,\mu_r)$
\end{center} with $l>k$ and $\mu_l>0$. Subsequently, we suppose that \begin{center}
$v(D_k^iy_0)=(0,\ldots,0,\beta_{k_i},\ldots,\beta_r)$
\end{center} with $k_i>k$ and $\beta_{k_i}>0$. But \begin{center}
$D_k^{i+1}y_0=D_k(D_k^iy_0)=(d_{k_i}/d_k)D_{k_i}(D_k^iy_0)\sim \beta_{k_i}D_k^iy_0d_{k_i}/d_k$.
\end{center} So $v(D_k^{i+1}y_0)=v(D_k^{i}y_0)+v(d_{k_i}/d_k)$ and from Lemma \ref{prop:val_d_k}, $v(d_{k_i}/d_k)>\underline{0}$ with the expected property.

For the case $k>l$, we show that $i\in\{0,\ldots,n\},\ D_k^iy_0\sim\mu_l^iy_0(d_{l}/d_k)^i$. We suppose that $D_k^iy_0\sim\mu_l^iy_0(d_{l}/d_k)^i$ for some $i\in\{0,\ldots,n\}$. So  \begin{center}
$v(D_k^iy_0)=(0,\ldots,0,\mu_l,\mu_{l+1}+i(\theta_{l+1}^{(l)}- \theta_{l+1}^{(l)}), \ldots)$,
\end{center} which implies that  \begin{center}
$D_k^{i+1}y_0=D_k(D_k^iy_0)\sim D_k(\mu_l^iy_0(d_{l}/d_k)^i) \sim D_l(\mu_l^iy_0(d_{l}/d_k)^i)d_{l}/d_k\sim \mu_l^{i+1}y_0(d_{l}/d_k)^{i+1}$.
\end{center} So $v(D_k^iy_0)$ and $v(y_0)$ have the same sign.
\end{demo}

\subsection{Dealing with formal equations rather than polynomial ones.}\label{subsec:equ-formelle} 
There is an additional difficulty in dealing with a formal Differential Equation \ref{eq:equa-diff} rather than only a polynomial one: we have to verify when the evaluation $F(y_0,\ldots,y_0^{(n)})$ of the Differential Series \ref{eq:serie-diff} at some series $y_0$ is well-defined. 

We state without proof the following easy generalisation of a classical property \cite[Ch.VIII,Sec.5,Lemma]{fuchs:partial_ord}. 
\begin{defi}\label{defi:sommable}
Given an index set $I$, a family $\mathcal{F}=(a_i)_{i\in I}\in\mathds{K}_r^I$ is said to be 
\textsl{strongly summable} if:
\begin{itemize}
    \item $\textrm{Supp}\ \mathcal{F}:=\bigcup_{i\in I}\textrm{Supp}\ a_i$ is a well-ordered subset of  $\Gamma$;
\item for all $\alpha\in\textrm{Supp}\ \mathcal{F}$, the set 
$\{i\in I \ | \ \alpha\in \textrm{Supp}\ a_i\}$ is finite.
\end{itemize}
\end{defi}
\begin{lemma}
Given a strongly summable family $(a_i)_{i\in I}\in\mathds{K}_r^I$, then $\sum_{i\in I}a_i$ is well defined and, if we set $a_i=\sum_{\alpha\in\textrm{Supp}\ a_i}a_{i,\alpha} t^\alpha$ for all $i\in I$, then $ \sum_{i\in I}a_i=\sum_{\alpha\in\Gamma}(\sum_{i\in I}a_{i,\alpha})t^\alpha$.
\end{lemma} 

\begin{propo}\label{compo_cond1} Given a Differential Series
\ref{eq:serie-diff} and a generalised series $y_0\in\mathds{K}_r^\prec$ with $v(y_0^{(i)})>\underline{0}$ for all $i\in\{0,\ldots,n\}$, then $y_0$ is compatible with \ref{eq:serie-diff} (see Definition \ref{defi:stabil}).
\end{propo}
\begin{demo} We proceed by induction on $n$, considering an arbitrary $y_0\in\mathds{K}_r^\prec$ with $v(y_0^{(i)})>\underline{0}$ for all $i\in\{0,\ldots,n\}$. If $n=0$, we show that the family $(c_iy_0^i)_{i\in\mathbb{N}}$ is strongly summable. On the one hand, $\bigcup_{i\in\mathbb{N}}\textrm{Supp}\ (c_iy_0^i)\subset\bigcup_{i\in\mathbb{N}}\textrm{Supp}\ c_i+\left\langle \textrm{Supp}\ y_0\right\rangle$. By hypothesis $\bigcup_{i\in\mathbb{N}}\textrm{Supp}\ c_i=\textrm{Supp}\ F$ and $\textrm{Supp}\ y_0$ are well-ordered. Moreover $\textrm{Supp}\ y_0\subset\Gamma_{>0}$, so $\left\langle \textrm{Supp}\ y_0\right\rangle$ is also well-ordered. On the other hand, for any $\alpha\in\Gamma$ and $m\in\mathbb{N}$, the set of $(\gamma^{(m)},j_1\alpha_1+\cdots+j_l\alpha_l)\in\textrm{Supp}\ c_m\times \left\langle \textrm{Supp}\ y_0\right\rangle$ such that $j_1+\cdots+j_l=m$ and $\gamma^{(m)}+j_1\alpha_1+\cdots+j_l\alpha_l=\alpha$ is clearly finite.

If $n>0$, it suffices to consider a differential series $F(y,\ldots,y^{(n)})$ as an element of $\mathds{K}_r[[y,\ldots,y^{(n-1)}]][[y^{(n)}]]$ and then apply the induction hypothesis and the preceding lemma.
\end{demo}

\begin{coro}\label{coro:D_k-bien-definie}
Given a Differential series \ref{eq:serie-diff-D_k}, any series $y_0\in\mathds{K}_r^\prec$ is compatible with it.
\end{coro}

This last result is an immediate consequence of the Proposition \ref{deriv_posit}. The next one is follows from the Corollary \ref{coro:k_0}:
\begin{coro}\label{coro:der_log}
Given a generalised series $y_0\in\mathds{K}_r$, $v(y_0^{(i)})>\underline{0}$ for all $i\in\{0,\ldots,n\}$ if and only if:
\begin{equation}\label{eq:cond-posit}
v(y_0)>\max\{\underline{0},-nv(d_{k_0})\}.
\end{equation}
where $k_0$ is defined as in the Corollary \ref{coro:k_0}.
\end{coro}

\begin{remarque}\begin{itemize}
\item Note that in the case where  $v(d_{k_0})<\underline{0}$, if $y_0\in\mathds{K}_r$ is such that $v(y_0)>-nv(d_{k_0})$ as in the preceding corollary, then its analysis as in the Notation \ref{rem_decompo}, spells $y_0=y_{0,k}+\cdots+y_{0,r}$ with $k\geq k_0$.
\item This condition generalises the condition $\rho_i>n$ used in \cite{jcano:series} for the rational exponents $\rho_i$, $i\in\mathbb{N}$ of the series solution considered.
\end{itemize}
\end{remarque} 

\subsection{Controlling the supports.}
The purpose of this section is to understand the support of the evaluation of a  Differential Series \ref{eq:serie-diff} or \ref{eq:serie-diff-D_k} at some series $y_0=\sum_{\mu\in\textrm{Supp}\ y_0}m_\mu t^\mu\in\mathds{K}_r$. 
\begin{notation}\label{nota_initpart}
For any multi-index $I=(i_0,\ldots,i_n), J=(j_0,\ldots,j_n)\in\mathbb{N}^{n+1}$:\\
\noindent$\begin{array}{rcl}
\bullet\ \ |I|&=&i_0+\cdots+i_n\\
\|I\|&=&1i_1+2i_2+\cdots+ni_n\\
I!&=&i_0!\ i_1!\cdots i_n!;
\end{array}$\\
$\bullet$ $I+J$, $I- J$, denote respectively the termwise addition, substraction.\\
$\bullet$ We will use the classical partial ordering on $\mathbb{N}^{n+1}$: \begin{center}
$I\leq J \Leftrightarrow i_l\leq j_l$ for any $l=0,\ldots,n$.
\end{center}
$\bullet$  $F^{(I)}=\partial^{|I|}F/[\partial y^{i_0}\cdots(\partial y^{(n)})^{i_n}]$; \\
$\bullet$  For any initial segment $S$ of $\textrm{Supp}\ y_0$, we denote by $p_S$ the \textsl{initial part of $y_0$ with support $S$}, i.e. $p_S=\sum_{\mu\in S}m_\mu t^\mu$.\\
$\bullet$ For any $k=0,\ldots,r$, $p_S^{(I)_k}=p_S^{i_0}D_k(p_S)^{i_1}\cdots D_k^n(p_S)^{i_n}$. (Recall that $D_0$ is the original Hardy type derivation.) If the derivation is known from the context, we will denote simply $p_S^{(I)_k}=p_S^{(I)}$;\\
$\bullet$ $f_S^{(I)}=F^{(I)}(p_S)$ and $f_S=F(p_S)$;\\
$\bullet$ $v_S^{(I)}=v(F^{(I)}(p_S))$ and $v_S=v(F(p_S))$;\\
$\bullet$ For any proper initial segment $S\subsetneq\textrm{Supp}\ y_0$ and any successor segment $\tilde{S}$ of $S$ in $\textrm{Supp}\ y_0$ (\textsl{that is to say any initial segment of $\textrm{Supp}\ y_0\backslash S$}), we have by the \emph{Taylor expansion formula}:
\begin{equation}\label{taylor}
f_{S\bigcup\tilde{S}}=\sum_{I\in\mathbb{N}^{n+1}} \displaystyle\frac{f_S^{(I)}}{I!}p_{\tilde{S}}^{(I)}.
\end{equation}
\end{notation}
Now we introduce the following well-ordered subsets of $\Gamma_{>0}$ useful for the description of the supports of the differential series.
\begin{defi}\label{T_k}
For all $k\in\{1\ldots,r\}$, we define: \begin{center}
$\mathcal{T}_k=\displaystyle\sum_{i=1}^n\displaystyle\sum_{l=k}^r \left\langle\textrm{Supp}\ D_k^it_l/t_l\right\rangle$.
\end{center}
\end{defi}
\begin{lemma}\label{lemme:T_k}
For any $k\in\{1\ldots,r\}$, $\mathcal{T}_k$ is a well-ordered subset of  $\Gamma_{>0}$, obtained from $\textrm{Supp}\ t_k'/t_k$, $k=1,\ldots,r$, by finitely many elementary transformations.
\end{lemma}
\begin{demo}
It follows from Definition \ref{hypo_deriv} and the proof of Proposition \ref{deriv_posit}, that for all $k\geq l\in\{1,\ldots,r\}$ and $i\in\{1,\ldots,n\}$, $v(D_k^it_l/t_l)>\underline{0}$. So the $\mathcal{T}_k$'s are well-ordered subsets of $\Gamma_{>0}$.

From Remark \ref{rem:Hardy_deriv}, we deduce that, for any series $a\in\mathds{K}_r$ and any  $k=1,\ldots,r$:\begin{center}
$\begin{array}{lcl}
D_k(a)&=&a'/d_k\\
&=&[\sum_{\alpha} (a_\alpha\alpha_1) \alpha] t_1'/(t_1d_k)+\cdots+[\sum_{\alpha} (a_\alpha\alpha_r) \alpha] t_r'/(t_rd_k).
\end{array}$ 
\end{center}
Thus, $\textrm{Supp}\ D_k(a)\subset (\sum_{l=1}^r \textrm{Supp}\ a +\textrm{Supp}\ t_l'/t_l)-d_k$, and therefore is obtained from $\textrm{Supp}\ a$ and $\textrm{Supp}\ t_l'/t_l$, $l=1,\ldots,r$, by finitely many elementary transformations.
\end{demo}

\begin{propo}\label{supp_eval}
\begin{enumerate}
    \item Let $k\geq l\in\{1,\ldots,r\}$, a Differential Series \ref{eq:serie-diff-D_k} for the derivation $D_l$, and $y_{0,k}\in\mathds{K}_{r,k}$ be given. Then:\begin{center}
 $\textrm{Supp}\
(F(y_{0,k},\ldots,D_l^ny_{0,k}))\subset \textrm{Supp}\
F+\mathcal{T}_l+\left\langle\textrm{Supp}\ y_{0,k}\right\rangle$.
\end{center}
\item Consider a Differential Series \ref{eq:serie-diff}, and a compatible series $y_{0,k}\in\mathds{K}_{r,k}$ for some $k\in\{1,\ldots,r\}$. Then, with the notations of Corollary \ref{coro:k_0}:
\begin{enumerate}
\item if $k< k_0$, then we have:\begin{center}
$\textrm{Supp}\
(F(y_{0,k},\ldots,y_{0,k}^{(n)}))\subset \textrm{Supp}\
F+\mathcal{T}_k+\left\langle\bigcup_{i=0}^n\textrm{Supp}\ y_{0,k}+iv(d_k)\right\rangle$;\end{center}
\item if $k\geq k_0$, then we have:\begin{center}
$\textrm{Supp}\
(F(y_{0,k},\ldots,y_{0,k}^{(n)}))\subset \textrm{Supp}\
F+\mathcal{T}_{k_0}+\left\langle\bigcup_{i=0}^n\textrm{Supp}\ y_{0,k}+iv(d_{k_0})\right\rangle$\end{center}
and 
\begin{center}
 $\textrm{Supp}\
(F(y_{0,k},\ldots,y_{0,k}^{(n)}))\subset (\bigcup_{I}\textrm{Supp}\ c_I+\|I\|v(d_{k_0})+\mathcal{T}_{k_0}+\left\langle\textrm{Supp}\ y_{0,k}\right\rangle$.
\end{center}
\end{enumerate}
\end{enumerate}
\end{propo}
\begin{demo} We treat the two cases at a time, by taking $l\in\{0,\ldots,r\}$. 
For any $I\in\mathbb{N}^{n+1}$, we have $\textrm{Supp}\ (c_Iy_{0,k}^{(I)})\subset \textrm{Supp}\
c_I+\textrm{Supp}\ y_{0,k}^{i_0}+\ldots+\textrm{Supp}\ (D_l^ny_{0,k})^{i_n}$. But, for any element $\mu=(0,\ldots,0,\mu_k,\ldots,\mu_r)$ with $\mu_k\neq 0$ of $\Gamma$, we have: \begin{center}
$D_l(t^\mu)=t^\mu(\mu_kD_lt_k/t_k\ +\cdots+\mu_rD_lt_r/t_r)$.
\end{center}
By induction, one can deduce the following intermediate result (for a detailed proof, see \cite[Lemma 4.2.12]{matu:these}):
\begin{lemma}\label{lemma_eval1}
Let $k\in\{1,\ldots,r\}$, $l\in\{0,\ldots,r\}$, and
$i\in\{1,\ldots,n\}$. Let $y_{0,k}\in\mathds{K}_{r,k}$. Then\begin{center}
$\textrm{Supp}\ D_l^iy_{0,k}\subset \textrm{Supp}\ y_{0,k}+
\bigcup\textrm{Supp}\
[(D_lt_k/
t_k)^{j^{(k)}_1}\cdots(D_l^it_k/
t_k)^{j^{(k)}_i}\cdots (D_lt_r/
t_r)^{j^{(r)}_1}\cdots(D_l^it_r/
t_r)^{j^{(r)}_i}]$ 
\end{center}
 where the union is taken over
\begin{center}
 $j^{(k)}_1+\cdots+ij^{(k)}_i+\cdots+j^{(r)}_1+ \dots+ij^{(r)}_i=i$.
\end{center}
\end{lemma}

In the case where $k\geq l\in\{1,\ldots,r\}$, we deduce that: \begin{center}
$\textrm{Supp}\ D_l^iy_{0,k}\subset \textrm{Supp}\ y_{0,k}+\mathcal{T}_k$.
\end{center}
Now, we note that, for any series $a\in\mathds{K}$ and any $j\in\mathbb{N}^*$, $\textrm{Supp }a^j\subset \langle\textrm{Supp }a\rangle$. So, for any $I\in\mathbb{N}^{n+1}$:\begin{center}
$\textrm{Supp}\ (c_Iy_{0,k}^{(I)_l})\subset 
 \textrm{Supp}\ c_I+\langle\textrm{Supp}\ y_{0,k}\rangle+\mathcal{T}_k$.
\end{center}
which implies the case (1).

For the case (2), we show by induction on $i$ that:\begin{description}
    \item[Case (a)] $\forall i=1,\ldots,n$ $\forall j\geq k$, $\textrm{Supp }(t_j^{(i)})/t_j\ \subset \mathcal{T}_{k}+i\theta^{(k)}$;

\item[Case (b)] $\forall i=1,\ldots,n$ $\forall j\geq k\geq k_0$, $\textrm{Supp }(t_j^{(i)})/t_j\ \subset \mathcal{T}_{k_0}+i\theta^{(k_0)}$.
\end{description}

Indeed, for $i=1$, we have $t_j'=d_kD_k(t_j)=d_{k_0}D_{k_0}(t_j)$. Suppose that the property holds for some $i=1,\ldots,n-1$. Thus, in the case (a), $t_j^{(i)}= d_{k}^it_j\sum_{m}m$ for some series $\sum_{m}m$ with support in \begin{center}
$\displaystyle\sum_{s=1}^i\displaystyle\sum_{p=k}^r \left\langle\textrm{Supp}\ D_k^st_p/t_p\right\rangle$,
\end{center}which is a subset of $\mathcal{T}_{k}$ (see Definition \ref{T_k}).
Therefore:\begin{center}
$\begin{array}{lcl}
t_j^{(i+1)}&=& d_{k}D_{k}(t_j^{(i)})\\
&=&id_{k}^iD_{k}(d_{k})t_j\sum_{m}m+d_{k}^{i+1}D_{k}(t_j)\sum_{m}m +d_{k}^{i+1}t_j\sum_{m}D_{k}(m)\\
&=& id_{k}^{i+1} t_j [\theta^{(k)}_{\tilde{k}}D_{k}(t_{\tilde{k}})/t_{\tilde{k}}\ +\cdots+ \theta^{(k)}_{r}D_{k}(t_{r})/t_r]\sum_{m}m+\\ &&\ d_{k}^{i+1}t_jD_{k}(t_j)/t_j\sum_{m}m +d_{k}^{i+1}t_j\sum_{m}D_{k}(m).
\end{array}$
\end{center}
Then it suffices to remark that for any $k\in\{1,\ldots,r\}$, for any monomial $m$ with support in $\mathcal{T}_{k}$, $D_{k}(m)$ has also its support in 
\begin{center}
$\displaystyle\sum_{s=1}^{i+1}\displaystyle\sum_{p=k}^r \left\langle\textrm{Supp}\ D_k^st_p/t_p\right\rangle$,
\end{center}
also included in $\mathcal{T}_{k}$.

In the case (b), $t_j^{(i)}= d_{k_0}^it_j\sum_{m}m$ for some series $\sum_{m}m$ with  support in
\begin{center}
$\displaystyle\sum_{s=1}^{i}\displaystyle\sum_{p=k_0}^r \left\langle\textrm{Supp}\ D_{k_0}^st_p/t_p\right\rangle$,
\end{center}
which is a subset of $\mathcal{T}_{k_0}$. Therefore:\begin{center}
$\begin{array}{lcl}
t_j^{(i+1)}&=& d_{k_0}D_{k_0}(t_j^{(i)})\\
&=&id_{k_0}^iD_{k_0}(d_{k_0})t_j\sum_{m}m+d_{k_0}^{i+1}D_{k_0}(t_j)\sum_{m}m +d_{k_0}^{i+1}t_j\sum_{m}D_{k_0}(m)\\
&=& id_{k_0}^{i+1} t_j [\theta^{(k_0)}_{k_0}D_{k_0}(t_{k_0})/t_{k_0}\ +\cdots+ \theta^{(k_0)}_{r}D_{k_0}(t_{r})/t_r]\sum_{m}m +\\&&\ d_{k_0}^{i+1}t_jD_{k_0}(t_j)/t_j\sum_{m}m +d_{k_0}^{i+1}t_j\sum_{m}D_{k_0}(m).
\end{array}$
\end{center}
The conclusion follows as for the preceding case.

Now, from the Lemma \ref{lemma_eval1}, we deduce that for any $I\in\mathbb{N}^{n+1}$:
\begin{description}
    \item[Case (a)] $\textrm{Supp }(c_Iy_{0,k}^{(I)}) \subset \textrm{Supp}\
c_I+\|I\|\theta^{(k)}+\textrm{Supp}\ y_{0,k}^{|I|}+\mathcal{T}_{k} $;
\item[Case (b)] $\textrm{Supp }(c_Iy_{0,k}^{(I)}) \subset \textrm{Supp}\
c_I+\|I\|\theta^{(k_0)}+\textrm{Supp}\ y_{0,k}^{|I|}+\mathcal{T}_{k_0} $.
\end{description}
\end{demo}

Note that, in the case (2)(a) ($k<k_0$), since $\tilde{k}>k$ (see Corollary \ref{coro:k_0}),  $v(y_{0,k})+i\theta^{(k)}$ has same sign as $v(y_{0,k})$ for any $i\in\mathbb{N}$. The following result provides some criterions of compatibility.

\begin{coro}\label{coro:equ-compat}
Let a Differential Series \ref{eq:serie-diff} and a generalised series $y_{0,k}\in\mathds{K}_{r,k}$ be given. The series $y_{0,k}$ is compatible with \ref{eq:serie-diff} if:\begin{description}
    \item[either] $k<k_0$
\item[or] $k\geq k_0$ and $v(y_{0,k})> \max\{\underline{0},-n\theta^{(k_0)}\}$
\item[or] $k\geq k_0$, $\theta^{(k_0)}<\underline{0}$, $\underline{0}< v(y_{0,k})\leq -n\theta^{(k_0)}$, and for all but finitely many $I\in\mathbb{N}^{n+1}$, $c_I=d_{k_0}^{\|I\|}a_I$ for some strongly summable family $(a_I)_{I}$, $a_I\in\mathds{K}_r$.
\end{description} 
\end{coro}

\subsection{The Weierstrass order of an equation.}
The following notion is a key one in our proof of Theorem \ref{theo_princ} in Section \ref{sect:weier=w}. It plays a role comparable to the one of Newton degree in the Newton polygon method \cite[Section 2.3.5]{vdh:autom-asymp}, \cite[Section 8.3.3]{vdh:transs_diff_alg}. We will need also to control its evolution when applying the transformations of the equations described in the next section.
\begin{defi}\textbf{(Weierstrass order)}\label{weier_defi}
A Differential Series \ref{eq:serie-diff} has Weierstrass order $w\in\mathbb{N}$ if:
\begin{itemize}
    \item for any $I\in\mathbb{N}^{n+1}$, $v(c_I)\geq\underline{0}$;
\item there exists $I\in\mathbb{N}^{n+1}$ with $|I|=w$ and $v(c_I)=\underline{0}$;
\item  for any $I\in\mathbb{N}^{n+1}$ with $|I|<w$, $v(c_I)>\underline{0}$.
\end{itemize} 
\end{defi}

Given a Differential Series \ref{eq:serie-diff} or \ref{eq:serie-diff-D_k}, if we divide it by $t^{\min(\textrm{Supp}\ F)}$ ($\min(\textrm{Supp}\ F$ exists since the support of $F$ is well-ordered), we obtain a series with Weierstrass order equal to the minimum of $|I|$ for multi-indexes $I\in\mathbb{N}^{n+1}$ such that $v(c_I)=\min(\textrm{Supp}\ F)$. Then, note that the support of the factored series is equal to $\textrm{Supp}\ F - \min(\textrm{Supp}\ F)$, which is the application of an elementary transformation  to $\textrm{Supp}\ F$.

So, \emph{without loss of generality, we will suppose from now that the differential series we consider has such a Weierstrass order}.

\section{Transformations of differential series.}\label{sec:transfo-series-diff}
\subsection{Changes of variable: additive and multiplicative conjugations.}
In \cite{matu:these}, we used these transformations respectively  under the name ``shiftings" and ``blow-ups". Then we became aware of, and now we resume the terminology "additive and multiplicative conjugations" for the same kind of transformations, which was introduced by van der Hoeven in \cite[Sections 5.2.2 and 5.2.3]{vdh:autom-asymp} and in \cite[Sections 8.2.1 and 8.2.2]{vdh:transs_diff_alg}. They play a role in particular in the proof of his Theorem 12.2 cited before.

We resume these two transformations and apply 
\begin{defi}\label{defi:transl}
Given a Differential Series  \ref{eq:serie-diff} or \ref{eq:serie-diff-D_k}, and some series $a\in\mathds{K}_r$, we denote $S_a=\textrm{Supp}\ a$. Then we call \textrm{additive conjugation} by $a$ the change of variable $y=a+\tilde{y}$ and we denote $$F_{S_a}(\tilde{y},\ldots,D_k^n\tilde{y}) =F(a+\tilde{y},\ldots,D_k^na+D_k^n\tilde{y})\ \ (k\in\{0,\ldots, r\})$$ the differential series thus obtained.
\end{defi}

\begin{propo}\label{transl:propos} 
\begin{enumerate}
\item Let $k\in\{1,\ldots,r\}$, a Differential Series \ref{eq:serie-diff-D_k} for the derivation $D_k$, and $a\in\mathds{K}_{r}^\prec$ with $a=a_k+\cdots+a_1$, be given. The differential series $F_{S_a}$ derived from \ref{eq:serie-diff-D_k} by additive conjugation by $a$ has a well-ordered support. Moreover:\begin{center}
 $\textrm{Supp}\ F_{S_a}\subset \textrm{Supp}\
F+\mathcal{T}_k+\left\langle S_a\right\rangle$.
\end{center}
\item Let a Differential Series \ref{eq:serie-diff} and a series $a\in\mathds{K}_r^\prec$ with $v(a^{(i)})>\underline{0}$ for all $i=0,\ldots,n$, be given. The differential series $F_{S_a}$ derived from \ref{eq:serie-diff} by additive conjugation by $a$, has a well-ordered support. Moreover, if $a=a_l\in\mathds{K}_{r,l}$ for some $l=1,\ldots,r$, then we consider the following cases:\begin{enumerate}
\item if $l< k_0$, then we have:\begin{center}
$\textrm{Supp}\
F_{S_a}\subset \textrm{Supp}\
F+\mathcal{T}_l+\left\langle\bigcup_{i=0}^nS_a+iv(d_l)\right\rangle$;\end{center}
\item if $l\geq k_0$, then we have:\begin{center}
$\textrm{Supp}\
F_{S_a}\subset \textrm{Supp}\
F+\mathcal{T}_{k_0}+\left\langle\bigcup_{i=0}^nS_a+iv(d_{k_0})\right\rangle$\end{center}
and 
\begin{center}
 $\textrm{Supp}\
F_{S_a}\subset (\bigcup_{I}\textrm{Supp}\ c_I+\|I\|v(d_{k_0}))+\mathcal{T}_{k_0}+\left\langle S_a\right\rangle$.
\end{center}
\end{enumerate}
\end{enumerate}
\end{propo}
\begin{demo} By the additive conjugation $y=\tilde{y}+a$, we have in the two cases of the proposition $D_k^iy=D_k^i\tilde{y}+D_k^ia$ for any $i\in\{1,\ldots,n\}$, $k\in\{0,\ldots,r\}$. By the Taylor expansion formula \ref{taylor}, we have:
\begin{center}
$F_{S_a}(\tilde{y},\ldots,D_k^n\tilde{y})= F(\tilde{y}+a,\ldots,D_k^n\tilde{y}+D^na)=\sum_{J\in\mathbb{N}^{n+1}} 
(f^{(J)}_{S_a}/J!)\tilde{y}^{(J)}.$
\end{center}
Note that $F^{(J)}=\sum_{I\geq J} c_Iy^{(I-J)}$, which implies that
$\textrm{Supp}\ F^{(J)}\subset \textrm{Supp}\ F$ for any $J$. By hypothesis, $a$ is compatible with any differential series $f^{(J)}$,  $J\in\mathbb{N}^{n+1}$: $f^{(J)}_{S_a}=F^{(J)}(a,\ldots,D_k^na)$ is a well-defined element of $\mathds{K}_r$. Moreover, from Proposition \ref{supp_eval}, we deduce that for any $J\in\mathbb{N}^{n+1}$:
\begin{enumerate}
\item if $k\neq 0$, then $\textrm{Supp}\
f^{(J)}_{S_a}\subset \textrm{Supp}\
F+\mathcal{T}_k+\left\langle S_a\right\rangle$;
\item if $k=0$, then:
\begin{enumerate}
    \item if $l<k_0$, then: $\textrm{Supp}\
f^{(J)}_{S_a}\subset \textrm{Supp}\
F^{(J)}+\mathcal{T}_k+\left\langle\bigcup_{i=0}^nS_a+iv(d_k)\right\rangle$\\
\item if $l\geq k_0$, then:  $\textrm{Supp}\
f^{(J)}_{S_a}\subset \textrm{Supp}\
F^{(J)}+\mathcal{T}_{k_0}+\left\langle\bigcup_{i=0}^nS_a+iv(d_{k_0}) \right\rangle$\\ 
and also \\
$\textrm{Supp}\
f^{(J)}_{S_a}\subset(\bigcup_{c_I\neq 0,\ I\geq J}\textrm{Supp}\ c_I+\|I-J\|v(d_{k_0}))+\mathcal{T}_{k_0}+\left\langle S_a\right\rangle$.
\end{enumerate}
\end{enumerate}
So, we deduce for $\textrm{Supp}\ F_{S_a}=\cup_{J} \textrm{Supp}\
f^{(J)}_{S_a}$, the desired properties.
\end{demo}

\begin{remark}\label{rem:weier_transl}
Consider a Differential Series \ref{eq:serie-diff} or \ref{eq:serie-diff-D_k} and a series $y_0\in\mathds{K}_r^\prec$ such that $v(y_0^{(i)})>\underline{0}$. for any $i=0,\ldots,n$. If we denote $F=\sum_{I\in\mathbb{N}^{n+1}}c_Iy^{(I)}$ and $F_{S_a}=\sum_{I\in\mathbb{N}^{n+1}}c'_I\tilde{y}^{(I)}$, we have $c_I=F^{(I)}(0)/I!$ and $c'_I=f^{(I)}_{S_a}/I!$. So, $v(c_I)=\underline{0}$ if and only if $v(c'_I)=\underline{0}$, and for any such $I$'s, $\delta(c_I)=\delta(c'_I)$. Therefore,\emph{ the Differential Series \ref{eq:serie-diff} and $F_{S_a}$ have same Weierstrass order}.
\end{remark}

In the sequel, we will use \textsl{anti-lexicographical ordering} for multi-indexes: for any $I=(i_0,\ldots,i_n)$ and any $J=(j_0,\ldots,j_n)\in\mathbb{N}^{n+1}$,\begin{center}
 $I<_{antilex}J \Leftrightarrow \textrm{Supp}\ (I-J)\neq\emptyset$ and $i_k<j_k$ where $k=\max(\textrm{Supp}\ (I-J))$.
\end{center}

\begin{defi}\label{defi:eclt} Let $m=t^\lambda$, $\lambda\in\Gamma$, be a monic monomial. We call multiplicative conjugation by $m$ the change of variable $y=mz$.
\end{defi}
\begin{propo}\label{eclt:propos} 
\begin{enumerate}
\item Let  a Differential Series \ref{eq:serie-diff-D_k} and a monic monomial $m=t^\lambda$ with  $\lambda=(0,\ldots,0,\lambda_l,\ldots,\lambda_r)$, $\lambda_k>0$, for some $l=k,\ldots,r$, be given. Performing the multiplicative conjugation $y=mz$, we obtain a differential series \begin{center}
$\hat{F}(z,\ldots,D_k^nz)=\sum_{J\in\mathbb{N}^{n+1}} \hat{c}_Jz^{(J)_k}$
\end{center} with well-ordered support such that:
\begin{center}
$\textrm{Supp}\ \hat{F}\subset\textrm{Supp}\ F+\mathbb{N}\lambda+\mathcal{T}_k$.
\end{center}
\item Let a Differential Series \ref{eq:serie-diff} and a monic monomial  $m=t^\lambda$ with $\lambda=(0,\ldots,0,\lambda_l,\ldots,\lambda_r)$, $\lambda\neq 0$ be given. Performing the multiplicative conjugation $y=mz$, we obtain a differential series \begin{center}
$\hat{F}(z,\ldots,z^{(n)})=\sum_{J\in\mathbb{N}^{n+1}} \hat{c}_Jz^{(J)}$
\end{center} such that: 
\begin{enumerate}
\item if $l< k_0$, then we have:\begin{center}
$\textrm{Supp}\ \hat{F} \subset \textrm{Supp}\
F+\mathcal{T}_l+\sum_{i=0}^n\mathbb{N}(\lambda+iv(d_l))$;\end{center}
\item if $l\geq k_0$, then we have:\begin{center}
$\textrm{Supp}\
\hat{F}\subset \textrm{Supp}\
F+\mathcal{T}_{k_0}+\sum_{i=0}^n\mathbb{N}(\lambda+iv(d_{k_0}))$\end{center}
and 
\begin{center}
$\textrm{Supp}\ \hat{F}\subset (\bigcup_{c_I\neq 0}\bigcup_{0\leq j\leq \|I\|}\textrm{Supp}\ c_I+|I|\lambda+jv(d_{k_0}))+\mathcal{T}_{k_0}$.
\end{center}
\end{enumerate}
If $m$ is strongly compatible with \ref{eq:serie-diff}, then $\hat{F}$ has well-ordered support.
\end{enumerate}
\end{propo}
\begin{demo}
We treat the two cases at a time, taking $k=0,\ldots,r$. For any $j=0,\ldots,n$, $D_k^jy=\sum_{i=0}^jC_j^iD_k^{j-i}mD_k^iz$. So, for any $I\in\mathbb{N}^{n+1}$, \begin{center}
$y^{(I)}=(mz)^{i_0}(D_kmz+mD_kz)^{i_1}\cdots(\sum_{i=0}^n C_n^iD_k^{n-i}mD_k^iz)^{i_n}= \sum l_{J,K}m^{(K)}z^{(J)}$
\end{center} where $l_{J,K}\in\mathbb{N}^*$ and the sum is taken over $J,K\in\mathbb{N}^{n+1}$ such that $|J|=|K|=|I|$, $\|K\|+\|J\|=\|I\|$ and $J,K\leq_{antilex}I$. We deduce that:
\begin{lemma}\label{lemme:coeff_eclt}
 For any $J\in\mathbb{N}^{n+1}$, 
 $\hat{c}_J=\sum k_{I,K}c_Im^{(K)}$ where $k_{I,J}\in\mathbb{N}^*$ and the sum is taken over the $I,K\in\mathbb{N}^{n+1}$ such that 
$|I|=|J|=|K|$, $\|I\|=\|K\|+\|J\|$ and $J,K\leq_{antilex}I$.
\end{lemma}
Therefore, we have:
\begin{center}
$\textrm{Supp }\hat{c}_J=\bigcup_{I,J,K}\textrm{Supp }c_I + \textrm{Supp }m^{(K)}$.
\end{center}
Then the result is obtained as for the proof of the Proposition \ref{supp_eval}.
\end{demo}

\begin{remark}\label{rem:weier_eclt}
In the case where $l=k\in\{1,\ldots,r\}$ and $\lambda_k>0$, we have for any multi-index $K$, $m^{(K)}\sim\lambda_l^{\|K\|}m^{|K|}$. Then, by Lemma \ref{lemme:coeff_eclt}, we obtain that \begin{center}
$\hat{c}_J=\sum k_{I,J}c_I\lambda_l^{\|K\|}m^{|K|} =\sum k_{I,J}\lambda_l^{\|K\|}c_Im^{|I|}$. 
\end{center} It means that a term $D_k^iy$ generates by multiplicative conjugation an analogous term $mD_k^iz$, plus terms with order of derivation in $z$ lower than $i$. Fix $l\in\mathbb{N}$ and consider the terms $c_Iy^{(I)}$ of $F$ with $|I|=l$ (suppose that there exists at least one). We denote $v_0=\min\{v(c_I)\ |\ |I|=k\}$, and \begin{center}
$\mathcal{A}_l=\{I\in\mathbb{N}^{n+1}\ |\ v(c_I)=v_0,\ |I|=l\}$.
\end{center} Then the terms $c_Iy^{(I)}$ for $I\in\mathcal{A}_l$ provide by multiplicative conjugation at least one term $\hat{c}_{I_0}z^{(I_0)}$ with $v(\hat{c}_{I_0})=v_0+lv(M)$: the one for $I_0=\max_{anti-lex}(\mathcal{A}_l)$ (its leading term can not be cancelled by any other term).

In particular, for a differential series with Weierstrass order $w$, we consider $\mathcal{A}_w$ the set of multi-indexes $I$ for which $v(c_I)=\underline{0}$. Then there exists a coefficient $\hat{c}_{I_0}$ with $|I_0|=w$ and with valuation $wv(M)$. Since for any $I$ with $|I|>w$, we have $v(\hat{c}_I)\geq (w+1)v(M)$ which is bigger than $w.v(M)$, then, denoting \begin{center}
$v_{min}=\min\{v(\hat{c}_I)\ |\ I\in\mathbb{N}^{n+1}\}$,
\end{center} we have $v_{min}\leq w.v(M)$. Thus \emph{if we divide the new series $\hat{F}$ by $t^{v_{min}}$, we obtain a series with Weierstrass order $\hat{w}$ at most equal to $w$ the Weierstrass order of the initial series $F$}.
\end{remark}

\subsection{Changes of derivation.}\label{sect:chgement-deriv}
Given a Differential Series \ref{eq:serie-diff} or \ref{eq:serie-diff-D_k}, denoted equally with $D_k$ for some $k\in\{0,\ldots,r-1\}$, can we transform it into some $\tilde{F}_{k,l}(y,\ldots,D_l^ny)$ for some $l\in\{k+1,\ldots,r\}$ ?  To do this, on the one hand we need to express the transformation, i.e. find formal formulas connecting $F$ and $\tilde{F}_{k,l}$. This is the purpose of Proposition \ref{chang_deriv}. On the other hand, since it may happen that $\tilde{F}_{k,l}$ does not have a well-ordered support, we must check when such a transformation is well defined (Proposition \ref{chgt_cond1}). Finally we determine the support of such a well-defined $\tilde{F}_{k,l}$ (Proposition \ref{supp_chang-deriv0}).

\begin{propo}\label{chang_deriv}
Let $k,l\in\{0,\ldots,r\}$, we denote
$m=d_l/d_k$ (setting $d_0=1$). We set for any
$i\in\{1,\ldots,n\}$,\begin{center}
$D_k^iy=q_{1,i}D_ly+q_{2,i}D_l^2y+\ldots+q_{i,i}D_l^iy$
\end{center}  where the 
multi-sequence $(q_{j,i})_{j=1,..,i;i=1,..,n}$ is defined by
:\begin{center}
$\left\{ \begin{array}{lclll}
q_{1,1}&=&m,\\
q_{1,i+1}&=&D_kq_{1,i},& \ \forall i=1,\cdots, n-1, &\ \\
q_{j+1,i+1}&=&q_{j,i}m+D_kq_{j+1,i},& \ \forall i=1,\cdots, n-1,& \ \forall j=1,\cdots, i-1,\\
q_{i+1,i+1}&=&q_{i,i}m,& \ \forall i=1,\cdots, n-1, &\ \\
\end{array}\right.$
\end{center}
Then we denote
$F(y,\ldots,D_k^ny)
=\tilde{F}_{k,l}(y,\ldots,D^n_ly)$ (note that $\tilde{F}_{k,l}$ may not have a well-ordered support). 
In particular, $q_{1,i}=D_l^{i-1}m$ and $q_{i,i}=m^i$ for any $i\in\{1,\ldots,n\}$.
\end{propo}

\begin{demo}
We proceed by induction on $i$. If $i=1$, we have $D_ky=y'/d_k=(d_l/d_k). (y'/d_l)=(d_l/d_k)D_ly$. Thus we set $q_{1,1}=m$. For the induction, we suppose that \begin{center}
$D_k^iy=q_{1,i}D_ly+q_{2,i}D_l^2y+\ldots+q_{i,i}D_l^iy$.
\end{center} Therefore:\\
$\begin{array}{lcl}
D_k^{i+1}y&=&D_kq_{1,i}D_ly+q_{1,i}D_kD_ly+\ldots+D_kq_{i,i}D_l^iy+q_{i,i}D_kD_l^iy\\
\ &=&D_kq_{1,i}D_ly+q_{1,i}(d_l/
d_k)D_l^2y+\ldots+D_kq_{i,i}D_l^iy+q_{i,i}(d_l/d_k)D_l^{i+1}y\\
\ &=&D_kq_{1,i}D_ly+(q_{1,i}m+D_kq_{2,i})D_l^2y+\ldots+ (q_{i-1,i}m+D_kq_{i,i})D_l^iy\\
\ & &+q_{i,i}mD_l^{i+1}y.
\end{array}$\\
The property holds at the step $i+1$.
\end{demo}

\begin{lemma}\label{multi-suite} For any $i\in\{1,\ldots,n\}$ and $j\in\{1,\ldots,i\}$, $q_{j,i}$ is a differential polynomial in $m$ such that
$q_{j,i}=\sum q_{j,i,I}m^{(I)_k}
\ \in\mathds{K}_r$ where the sum is taken over $I\in\mathbb{N}^i$
with $|I|=j$ and $\|I\|=i-j$.
\end{lemma}
\begin{demo}
We proceed by induction on $i$. If $i=1$, $q_{1,1}=m$ and we set
$q^{1,1,(1)}=1$. If $i=2$, $q_{1,2}=D_{k}m$ and $q_{2,2}=m^2$. So  $q_{1,2,(0,1)}=q_{2,2,(2,0)}=1$.

For the induction, we suppose that $q_{1,i+1}=D_{k}^im, \
q_{1,i+1,(0,..,1)}=1$ and $q_{i+1,i+1}=m^{i+1}, \
q_{i+1,i+1,(i+1,0,..)}=1$. Then \begin{center}
$q_{j+1,i+1}=mq_{j,i}+D_{k}q_{j+1,i}=\sum_{|I|=j,\|I\|=i-j}
q_{j,i,I}m^{(I+(1,0,..,0))_k} +\sum_{|I|=j+1,\|I\|=i-j-1}
q_{j,i,I}D_{k}(m^{(I)_k}).$
\end{center}
In the right hand term, for the first sum, we set
$J=I+(1,0,..,0)$. So $|J|=j+1$ and
$\|J\|=i+1-j-1$. For the second one, we have \begin{center}
$D_{k}(m^{(I)_k})=D_{k}(m^{i_0}(D_{k}m)^{i_1} \ldots(D_{k}^{i-1}m)^{i_{i-1}})=\sum_{l=0}^{i-1}i_l m^{i_0}\ldots(D_{k}^lm)^{i_l-1}(D_{k}^{l+1}m)^{i_l+1}\ldots (D_{k}^{i-1}m)^{i_{i-1}}.$
\end{center}
We set $K=I+ (0,..,-1,1,..,0)$ with -1 in
$l^{th}$ position. So $|K|=j+1$ and $\|K\|=i+1-j-1$. Then we obtain, for any $l=1,\ldots,i-1$,
\begin{center}
$q_{j+1,i+1,I}=q_{j,i,I-(1,0,..,0)}+(i_l+1) q_{j,i,I-(0,..,-1,1,..,0)}$
\end{center}
and for all $I\in\mathbb{N}^{i+1}$ with $|I|=j+1$, $\|I\|=i+1-j-1$,
$i_1\geq 1$ and $i_{l+1}\geq 1$.
\end{demo}

As an example, we provide in the following table some polynomials $q_{j,i}$.\\
\begin{tabular}{|c|c|c|c|c|}
  \hline $i\backslash j$& 1 & 2 & 3 & 4 \\\hline
  1 &$m$ & - & - & -  \\\hline
  2 & $D_{k}m$ & $m^2$ & - & - \\\hline
  3 & $D_{k}^2m$ & $3mD_{k}m$ & $m^3$ & - \\\hline
  4 & $D_{k}^3m$ & $4mD_{k}^2m+3(D_{k}m)^2$ & $6D_{k}m.m^2$ & $m^4$ \\ \hline
  5 & $D_{k}^4m$ & $5mD_{k}^3m+10D_{k}mD_{k}^2m$ & $10m^2D_{k}^2m+15m(D_{k}m)^2$ & $10m^3D_{k}m$ \\ \hline
\end{tabular}

\begin{remark}\label{rem:chge_deriv}
From Proposition \ref{chang_deriv} and Lemma \ref{multi-suite}, $D_k^iz=\sum_{j=1}^iq_{i,j}D_l^jz$ with $q_{i,j}=\sum q_{j,i,I}m^{(I)_k}\in\mathds{K}_r^\prec$ where $m=d_l/d_k$, $q_{j,i,I}\in\mathbb{N}$ and the sum is taken over multi-indexes $I\in\mathbb{N}^i$ such that $|I|=j$ and $\|I\|=i-j$. In particular, the coefficient of $D_l^iz$ is $q_{i,i}=m^i$.
 So, given some differential monomials $a_Iy^{(I)_k}$ with coefficients that have same valuation $v^{(0)}$, by the change of derivation they generate new differential monomials $b_Jy^{(J)_l}$ with $|I|=|J|$ and with one of them (the one with $J=I_0$ where $I_0$ is the greatest multi-index for the anti-lexicographical ordering among these $I$'s) with coefficient that has valuation  $v^{(0)}+\|I_0\|v(m)$. But $v(m)=(0,\ldots,0,\mu_{l_1},\ldots,\mu_r)$ with $\mu_{l_1}>0$ and $l_1>k$ (see Proposition \ref{prop:val_d_k}). So such a coefficient has a valuation of type \begin{center}
$v^{(0)}+\|I\|(0,\ldots,0,\mu_{l_1},\ldots,\mu_r)$ 
\end{center}with $\mu_{l_1}>0$ and $l_1>k$. \emph{This remark will be useful to control the evolution of the Weierstrass order of the differential series through a change of derivation}.
\end{remark}

\begin{propo}\label{chgt_cond1}
\begin{enumerate}
\item We consider a Differential Series \ref{eq:serie-diff-D_k}. The changes of derivation \begin{center}
$F(y,\ldots,D_k^ny)=\tilde{F}_{k,l}(y,\ldots,D^n_ly)$
\end{center} are well defined for any $k<l$ in $\{1,\ldots,r\}$.
\item We consider a Differential Series \ref{eq:serie-diff}, the integer $k_0\in\{1,\ldots,r\}$ defined in the Corollary \ref{coro:k_0}, and some  $l\in\{1,\ldots,r\}$, with $l\leq k_0$ if $v(d_{k_0})<\underline{0}$ (see Corollary \ref{coro:der_log}).
\begin{enumerate}
\item If $v(d_l)\geq\underline{0}$, then, by the change of derivation  \begin{center}
$F(y,\ldots,y^{(n)})=\tilde{F}_{0,l}(y,\ldots,D^n_ly)$,
\end{center} we obtain a differential series $\tilde{F}_{0,l}$ with well-ordered support.
\item If $v(d_l)<\underline{0}$, then $l\leq k_0$. We apply successively to \ref{eq:serie-diff} the multiplicative conjugation $y=d_l^{-n}z$, and the change of derivation $F=\tilde{F}_{0,l}$. Then, $\tilde{F}_{0,l}$ has well-ordered support.
\end{enumerate}
\end{enumerate}
\end{propo}
\begin{demo}
Given $k,l\in\{0,\ldots,r\}$, by the Proposition \ref{chang_deriv}, if for any $i\in\{0,\ldots,n\}$, $v(D_k^im)\geq\underline{0}$ (where $m=d_l/d_k$), then by the change of derivation $F(y,\ldots,D_k^ny)=\tilde{F}_{k,l}(y,\ldots,D^n_ly)$, we obtain a differential series $\tilde{F}_{k,l}$ which has clearly a well-ordered support. But, in the case (1), $m= d_l/d_k$ with $l>k$. So $v(m)=v(d_l)-v(d_k)>\underline{0}$, and therefore, $v(D_k^im)>\underline{0}$ for all $i\in\{0,\ldots,n\}$ (see Proposition \ref{deriv_posit}). 

In the case (2)(a), $m=d_l$. By the Corollary \ref{coro:deriv_neg}, $v(m^{(i)})=v(d_l^{(i)})\geq\underline{0}$ for any $i$. 
 
For the case (2)(b), we denote $M=d_l^{-n}$. Combining the results of the Propositions \ref{eclt:propos} and \ref{chang_deriv}, we show that:
\begin{lemma}
For any $i\in\mathbb{N}$, $y^{(i)}$ becomes a linear combination with positive integer coefficients, of terms
\begin{center}
$M^{(i-k-j)}d_l^{(I)}D_l^kz= M^{(i-k-j)}d_l^{i_0}(d_l')^{i_1}\cdots(d_l^{(j)})^{i_j}D_l^kz$
\end{center} where $k=0,\ldots,i$, $j=0,\ldots,i-k$, $I\in\mathbb{N}^{j+1}$ with $|I|=k$ and $\|I\|=j$. 
\end{lemma}
\begin{demo}
Consider such a term $M^{(i-k-j)}d_l^{(I)}D_l^kz$. Then its derivative is:
\begin{center}
$\begin{array}{lcl}
(M^{(i-k-j)}d_l^{(I)}D_l^kz)'&=&M^{(i+1-k-j)}d_l^{(I)}D_l^kz+ M^{(i-k-j)}[i_0d_l^{(I+(-1,1,0,\cdots))}\\ &&+\cdots+i_jd_l^{(I+(0,\ldots,0,-1))}d_l^{(j+1)}]D_l^kz+\\ 
&&M^{(i-k-j)}d_l^{(I)}d_lD_l^{k+1}z\\
&=&M^{(i+1-k-j)}d_l^{(\tilde{I})}D_l^kz+ M^{(i+1-k-(j+1))}[i_0d_l^{(\tilde{I}+(-1,1,0,\cdots,0))}\\ &&+\cdots+i_jd_l^{(\tilde{I}+(0,\ldots,0,-1,1))}]D_l^kz\\ &&+ 
M^{(i+1-(k+1)-j)}d_l^{(\tilde{I}+(1,0,\ldots,0))}D_l^{k+1}z,\\
\end{array}$
\end{center}
where $\tilde{I}\in\mathbb{N}^{j+2} $ is defined as the multi-index $I$ expanded by the addition of a new component 0 at its end. All the new exponents have the desired properties.
\end{demo}

Recall that $M=d_l^{-n}$.  We show that:
\begin{lemma}
For any $m\in\mathbb{N}^*$, $M^{(m)}$ is an integral linear combination of terms 
\begin{center}
$d_l^{-n-p}d_l^{(Q)}=d_l^{-n-p}(d_l')^{q_1}\cdots (d_l^{(m-p)})^{q_{m-p+1}}$
\end{center}
where $p=1,\ldots,m$, $Q\in\mathbb{N}^{m}$ with $|Q|=p$ and $\|Q\|=m$.
\end{lemma}
\begin{demo}
By induction on $m$, we consider such a general term
$d_l^{-n-p}d_l^{(Q)}$. Then its derivative is:
\begin{center}
$\begin{array}{lcl}
(d_l^{-n-p}d_l^{(Q)})'&=&(-n-p)d_l^{-n-(p+1)}d_l^{(Q)}+ d_l^{-n-p}[q_1d_l^{(Q+(-1,1,0,\ldots))}+\cdots+\\ &&q_{m-p+1} d_l^{(Q+(0,\ldots,-1))}d_l^{(m-p+2)}]\\
&=&(-n-p)d_l^{-n-(p+1)}d_l^{(\tilde{Q})}+ d_l^{-n-p}[q_1d_l^{(\tilde{Q}+(-1,1,0,\ldots,0))}+\cdots+\\&& q_{m-p+1} d_l^{(\tilde{Q}+(0,\ldots,-1,1))}].
\end{array}$
\end{center}
where $\tilde{Q}\in\mathbb{N}^{m-p+2}$ is defined as the multi-index $Q$ expanded by the addition of a new component 0 at its end.  All the new exponents have the desired properties.
\end{demo}

 Now, combining the results of the two preceding lemmas, as the general term generated by  $y^{(i)}$ we obtain
\begin{center}
$M^{(i-k-j)}d_l^{(I)}D_l^kz= d_l^{-n-i+k+j}d_l^{i_0}(d_l')^{i_1+q_1}\cdots(d_l^{(m)})^{i_m+q_m}D_l^kz$
\end{center} 
where $m=\max\{j,i-k-j-p\}$ and the multi-indexes $I$ or $Q$ are completed by zeros when needed. But, by the Corollary \ref{coro:deriv_neg}, since $l\leq k_0$, $v(d_l^{(q)})\geq  (q+1)v(d_l)$ for any $q\in\mathbb{N}$. Thus we deduce that the valuation of the coefficient of this general term is bigger than the minimum of:
\begin{center}
$\begin{array}{lcl}
v(M^{(i-k-j)}d_l^{(I)})&\geq& v(d_l^{-n-i+k+j}d_l^{|I|+\|I\|+|Q|+\|Q\|})\\
&\geq &(-n+k+j+p)v(d_l)\\
&\geq &(-n+i) v(d_l)\\
&\geq &\underline{0} 
\end{array}$
\end{center}
since $|I|=k$, $\|I\|=j$, $|Q|=p$, $\|Q\|=i-k-j$, with $p=1,\ldots,i-k-j$, $j=0,\ldots,i-k$, $k=0,\ldots,i$ and $i=0,\ldots,n$.
\end{demo}

\begin{propo}\label{supp_chang-deriv0}
For any  $k<l\in\{0,\ldots,r\}$, with the same hypothesis as in the Proposition \ref{chgt_cond1}, we have:\begin{enumerate}
\item $\textrm{Supp}\ \tilde{F}_{k,l} \subset\textrm{Supp}\ F+\mathcal{T}_k$
\item \begin{enumerate}
\item \begin{enumerate}
\item if $l<k_0$, then  $\textrm{Supp}\ \tilde{F}_{0,l}\subset\textrm{Supp}\ F+\mathbb{N}\ v(d_l)+\mathcal{T}_l$;
\item if $l\geq k_0$, then  $\textrm{Supp}\ \tilde{F}_{0,l}\subset\textrm{Supp}\ F+\mathbb{N}\ v(d_{k_0})+\mathcal{T}_{k_0}$;
\end{enumerate}
\item $\textrm{Supp}\ \tilde{F}_{0,l}\subset\textrm{Supp}\ F-\mathbb{N}v(d_l)+\mathcal{T}_l$.
\end{enumerate}
\end{enumerate}
\end{propo}
\begin{demo}
(1) We denote $m=d_l/d_k$. From Lemma \ref{prop:val_d_k}, we have \begin{center}
$\textrm{Supp}\ m=\{v(m)\}=\{v(d_l)-v(d_k)\}=\{(0,\ldots,0,\beta_j,\ldots,\beta_r)\}$
\end{center} for some $j>k$ and some $\beta_j,\ldots,\beta_r\in\mathbb{R}$ with $\beta_j>0$.
The change of derivation reads $D_k^iy=\sum_{j=1}^iq_{i,j}D_l^jy$ for any $i\in\{0,\ldots,n\}$, with \begin{center}
$q_{i,j}=\sum_{I\in\mathbb{N}^i,|I|=j,\|I\|=i-j}q^{j,i,I} m^{(I)_k}$.
\end{center} From Lemma \ref{supp_eval}, for any $i\in\{0,\ldots,n\}$, we have $\textrm{Supp}\ D_k^im\subset v(m)+\mathcal{T}_k$. Then it suffices to note that $m=d_l/d_k\sim t'_l/t_ld_k= D_kt_l/t_l$ and so $\textrm{Supp}\ m\subset\mathcal{T}_k$.

(2)(a) If $v(d_l)\geq\underline{0}$, by the Propositions \ref{supp_eval}, \ref{chang_deriv} and the Lemma \ref{multi-suite}, we have, for any $i\in\{1,\ldots,n\}$ and any $j\in\{1,\ldots,i\}$:
 \begin{center}
$\begin{array}{rllcl}
\mathrm{(i)}& \mathrm{if\ l<k_0,\ then\ }&\textrm{Supp}\ q_{i,j}&\subset&\textrm{Supp}\ d_l+\mathbb{N}v(d_l)+\mathcal{T}_l;\\
\mathrm{(ii)}& \mathrm{if\ l\geq k_0,\ then\ }&\textrm{Supp}\ q_{i,j}&\subset&\textrm{Supp}\ d_l+\mathbb{N}v(d_{k_0})+\mathcal{T}_{k_0}.
\end{array}$
\end{center}
We have $\textrm{Supp}\ d_l=\{v(d_l)\}\subset \textrm{Supp}\  t_l'/t_l$. But,  $t_l'/t_l=d_{k_0}D_{k_0}t_l/t_l$. So, $v(d_l)\in v(d_{k_0})+\mathcal{T}_{k_0}$. We deduce that:
 \begin{center}
$\begin{array}{rllcl}
\mathrm{(i)}& \mathrm{if\ l<k_0,\ then\ }&\textrm{Supp}\ q_{i,j}&\subset&\mathbb{N}v(d_l)+\mathcal{T}_l;\\
\mathrm{(ii)}& \mathrm{if\ l\geq k_0,\ then\ }&\textrm{Supp}\ q_{i,j}&\subset& \mathbb{N}v(d_{k_0})+\mathcal{T}_{k_0},
\end{array}$
\end{center}
which leads directly to the desired conclusion.

(2)(b) As in the proof of the Proposition \ref{chgt_cond1}, by combination of the multiplicative conjugation $y=d_l^{-n}$ and the change of derivation $F=\tilde{F}_{0,l}$, a term $y^{(i)}$ for some $i=0,\ldots,n$ generates an integral linear combination of terms \begin{center}
$d_l^{-n-i+k+j}d_l^{i_0}(d_l')^{i_1+q_1}\cdots(d_l^{(m)})^{i_m+q_m}D_l^kz$
\end{center} 
where $|I|=k$, $\|I\|=j$, $|Q|=p$, $\|Q\|=i-k-j$, with $p=0,\ldots,i-k-j$, $j=0,\ldots,i-k$, $k=0,\ldots,i$.

We rewrite such a term as:
\begin{center}
$d_l^{-n-i+k+j}d_l^{|I|+\|I\|+|Q|+\|Q\|} (d_l'/d_l^2)^{i_1+q_1}\cdots(d_l^{(m)}/d_l^{m+1})^{i_m+q_m}D_l^kz= d_l^{-n+p+k+j}(d_l'/d_l^2)^{i_1+q_1}\cdots(d_l^{(m)}/d_l^{m+1})^{i_m+q_m}D_l^kz.$
\end{center}
Note that $0\leq p+k+j\leq n$.

Now, it suffices to recall that show that:
\begin{lemma}
For any $s\in\mathbb{N}^*$, $\textrm{Supp }d_l^{(s)}\subset (s+1)v(d_l)+\mathcal{T}_l$.
\end{lemma}
\begin{demo}
We proceed by induction on $s$. Since $l\leq k_0$, we have for some $\tilde{l}\geq l$, $d_l'=d_l(\theta^{(l)}_{\tilde{l}} t_{\tilde{l}}'/t_{\tilde{l}}+\cdots+ \theta^{(l)}_rt_r'/t_r)$. So $d_l'=d_l^2(\theta^{(l)}_{\tilde{l}} D_l(t_{\tilde{l}})/t_{\tilde{l}}+\cdots+ \theta^{(l)}_rD_l(t_r)/t_r)$ which has a support as desired.

Suppose that, for some $s\in\mathbb{N}^*$, we have $d_l^{(s)}=d_l^{s+1}\sum c_\gamma t^\gamma$ for some series $\sum c_\gamma t^\gamma$ with support in $\mathcal{T}_l$. Then:
\begin{center}
$\begin{array}{lcl}
d_l^{(s+1)}&=&d_l^{s+1}[(s+1)(d_l'/d_l^1)\sum c_\gamma t^\gamma+\sum c_\gamma (t^\gamma)']\\
&=&d_l^{s+2}[(s+1)(d_l'/d_l^2)\sum c_\gamma t^\gamma+\sum c_\gamma D_l(t^\gamma)].\\
\end{array}$
\end{center}
\end{demo}
To conclude, note that $d_l'/d_l^2$ and $\sum c_\gamma D_l(t^\gamma)$ have support included in $\mathcal{T}_l$ (the first one by the result for $s=1$ here above, the second one by the case (1) of the Proposition \ref{supp_eval}). \end{demo}

\subsection{Reducing to positive valuation solutions for polynomial equations.}\label{sect:posit-val}
We consider a differential polynomial 
\begin{equation}\label{eq:serie-diff-poly}
P(y,\ldots,y^{(n)})=\sum_{I\in S\subset\mathbb{N}^{n+1},\ S\mathrm{\ finite}}
c_Iy^{(I)}\in\mathds{K}_r[y,\ldots,y^{(n)}]\setminus\{0\}.
\end{equation} We remind that in this case, any generalised series $y_0\in\mathds{K}_r$ is compatible with \ref{eq:serie-diff-poly} (see Definition \ref{defi:stabil}: the family $(c_Iy_0^{(I)})_I$, being finite, is strongly summable).
The purpose of this section is to derive from Theorem \ref{theo_princ} the same result, but for\textsl{ any} generalised series $y_0\in\mathds{K}_r$. 

We consider a generalised series \begin{center}
$y_0=m_0t^{\mu_0}+m_1t^{\mu_1}+\cdots\in\mathds{K}_r$
\end{center} where  $\delta(y_0)=m_0t^{\mu_0}$, such that $v(y_0^{(i)})\leq\underline{0}$ for some $i=0,\ldots,n$, i.e.  (see Corollary \ref{coro:der_log}) such that:
\begin{center}
$\mu_0\leq \alpha_0:=\max\{\underline{0},-nv(d_{k_0})\}$.
\end{center} 

\begin{lemma}\label{lemme:reduc-posit-val}
\begin{enumerate}
\item By the additive conjugation $y=\tilde{y}+m_0t^{\mu_0}$  (see Definition \ref{defi:transl}), we obtain from \ref{eq:serie-diff-poly} a new differential polynomial\begin{center}
$P_{\{\mu_0\}}(\tilde{y},\ldots,\tilde{y}^{(n)})=\sum_{I\in S}\tilde{c}_I\tilde{y}^{(I)}\in \mathds{K}_r[\tilde{y},\ldots,\tilde{y}^{(n)}] \setminus\{0\}$
\end{center}
with associated series \begin{center}
$\tilde{y}_0=y_0-m_0t^{\mu_0}=m_1t^{\mu_1}+\cdots\in\mathds{K}_r$.
\end{center}
\item By the multiplicative conjugation $\tilde{y}=t^{\mu_0-\alpha_0}\hat{y}$ (see Definition \ref{defi:eclt}), we obtain from $P_{\{\mu_0\}}$ a new differential polynomial \begin{center}
$\hat{P}(\hat{y},\ldots,\hat{y}^{(n)})=\sum_{J\in \hat{S} \ \mathrm{finite}}\hat{c}_J\hat{y}^{(J)}\in \mathds{K}_r[\hat{y},\ldots,\hat{y}^{(n)}] \setminus\{0\}$
\end{center}
with associated series \begin{center}
$\hat{y}_0=\tilde{y}_0/t^{\mu_0-\alpha_0} =m_1t^{\mu_1-\mu_0+\alpha_0}+\cdots\in\mathds{K}_r$.
\end{center}
which has valuation $v(\hat{y}_0)=\mu_1-\mu_0+\alpha_0>\alpha_0$.
\end{enumerate}
\end{lemma}
\begin{demo}
(1) By the Taylor expansion formula \ref{taylor}, we have \begin{center}
$\begin{array}{lcl}
P_{\{\mu_0\}}(\tilde{y},\ldots,\tilde{y}^{(n)})&=& P(y+m_0t^{\mu_0},\ldots,y^{(n)}+m_0(t^{\mu_0})^{(n)})\\
&=&\sum_{I\in S} \left(P^{(I)}(m_0t^{\mu_0},\ldots,m_0(t^{\mu_0})^{(n)})/I!\right) \tilde{y}^{(I)}.
\end{array}$
\end{center}
Since $P^{(I)}$ is a differential polynomial  for any $I$, we have $P^{(I)}(m_0t^{\mu_0},\ldots,m_0(t^{\mu_0})^{(n)})\in\mathds{K}_r$. Then $P_{\{\mu_0\}}$ is well-defined and has well-ordered support (which is the finite union of the supports of $\tilde{c}_I$).

\noindent(2) By the Lemma \ref{lemme:coeff_eclt}, we have $\hat{c}_J=\sum k_{I,K}\tilde{c}_I(t^{\mu_0-\alpha_0})^{(K)}$ where $k_{I,J}\in\mathbb{N}^*$ and the sum is taken over the $I,K\in\mathbb{N}^{n+1}$ such that 
$|I|=|J|=|K|$, $\|I\|=\|K\|+\|J\|$ and $J,K\leq_{antilex}I$. Therefore, we have a finite number of coefficients $\hat{c}_J$, any of which is itself a finite sum of generalised series $\tilde{c}_I(t^{\mu_0-\alpha_0})^{(K)}$. Thus, $\hat{P}$ is a differential polynomial.
\end{demo}

Thanks to the preceding lemma, in the case of a Differential Polynomial together with an arbitrary generalised series $y_0\in\mathds{K}_r$, we can reduce to the hypothesis of the Theorem \ref{theo_princ}. So, assuming that the latter holds, we obtain:

\begin{theo}\label{theo:poly-diff}
We denote $\alpha_0:=\max\{\underline{0},-nv(d_{k_0})\}$. Given a Differential Polynomial \ref{eq:serie-diff-poly} and a non zero series $y_0\in\mathds{K}_r$, there exists a well-ordered subset $\mathcal{R}$ of $\Gamma_{\geq 0}$ obtained from $\textrm{Supp}\ F$ and $\textrm{Supp}\ t_k'/t_k$, $k=1,\ldots,r$, by a finite number of elementary transformations such that:
\begin{itemize}
    \item either $\textrm{Supp}\ y_0\subseteq v(y_0)-\alpha_0+\mathcal{R}$;
\item or the series $y_0$ stabilises on $F$ with initial part $p_0$ and $\textrm{Supp}\ p_0\subseteq v(y_0)-\alpha_0+\mathcal{R}$. In this case, $y_0$ can not be a solution of the corresponding Differential Equation \ref{eq:equa-diff}.
\end{itemize}
\end{theo}
As a direct consequence, we obtain that:

\begin{theo}\label{theo:poly-diff1}
Given a differential equation \begin{equation}\label{eq:equa-diff-poly}
P(y,\ldots,y^{(n)})=0
\end{equation} where $P$ is a Differential Polynomial \ref{eq:serie-diff-poly}, the subset $(\textrm{Supp}\ y_0)-v(y_0)+\alpha_0$ of $\Gamma_{>\underline{0}}$, is obtained by finitely many elementary transformations from the supports of $F$ and the $t_k'/t_k$'s, $k=1,\ldots,r$. 
\end{theo}

In the case of a general Differential Series \ref{eq:serie-diff}, the changes of variable defined in the Lemma \ref{lemme:reduc-posit-val} may not generate a differential series $\hat{F}$ with well-ordered support. Nevertheless, we conjecture that the Theorems \ref{theo:poly-diff} and \ref{theo:poly-diff1} (or slighlty adapted versions of them) hold in this more general context.

 We are already working on this question, which we believe, will deserve a new paper by itself. 

\section{Proof of the Theorem \ref{theo_princ}}\label{sect:weier=w}
\subsection{The main lemma.}
The Theorem \ref{theo_princ} is a consequence of the following lemma:
\begin{lemma}\label{lemme:princ}
Given $k\in\{1,\ldots,r\}$ and $w\in\mathbb{N}$, we consider a Differential Series \ref{eq:serie-diff-D_k} and a generalised series $y_0=y_{0,k}+\cdots+y_{0,1}\in\mathds{K}_r^\prec$ with $y_{0,l}\in\mathds{K}_{r,l}$ for any $l\in\{1,\ldots,k\}$. Then there exists a well-ordered subset $\mathcal{R}$ of $\Gamma_{>0}$ obtained from $\textrm{Supp}\ F,\ \mathcal{T}_k,\ldots,\mathcal{T}_r$ by a finite number of elementary transformations such that:
\begin{itemize}
    \item either the exponents of $y_{0,k}$ belong to $\mathcal{R}$;
\item or the series $y_0$ stabilises on $F$ with initial part $p_0$ which is also a proper initial part of $y_{0,k}$. Moreover the exponents of $p_0$ belong to $\mathcal{R}$.
\end{itemize} 
\end{lemma}

\noindent\textit{Proof of Theorem \ref{theo_princ}.} Indeed, we consider a series $y_0\in\mathds{K}_r^\prec$ with  $v(y_0^{(i)})>\underline{0}$ for any $i=0,\ldots,n$, as in the statement of Theorem \ref{theo_princ}. We denote by 
$y_0=y_{0,k}+\cdots+y_{0,1}$ its decomposition as in the Notation \ref{rem_decompo}. Then, applying the point (2) of the Proposition \ref{chgt_cond1}, we consider two cases. If $v(d_k)\geq\underline{0}$, then, by the change of derivation  \begin{center}
$F(y,\ldots,y^{(n)})=\tilde{F}_{0,k}(y,\ldots,D^n_ky)$,
\end{center} the new differential series $\tilde{F}_{0,k}$ has well-ordered support. So it has  some Weierstrass order $w\in\mathbb{N}$ (see Definition \ref{weier_defi}). Moreover, by the Proposition \ref{supp_chang-deriv0}, the support of the new series $\tilde{F}_{0,k}$ is derived from the one of $F$ and the set $\mathcal{T}_k$ by elementary transformations. Thus, considering this series $\tilde{F}_{0,k}$ with $y_0$, we are reduced to the hypothesis of the Lemma \ref{lemme:princ}.

If $v(d_k)<\underline{0}$, as in the case (b) of the Proposition \ref{chgt_cond1}, we have $k\leq k_0$, and we can apply successively to \ref{eq:serie-diff} the multiplicative conjugation $y=d_k^{-n}z$, and the change of derivation $F=\tilde{F}_{0,k}$. Then, we obtain a  differential series $\tilde{F}_{0,k}$ with  well-ordered support. As above it has some Weierstrass order. Its support is also derived from $\textrm{Supp }F$ and $\mathcal{T}_k$ by elementary transformations. The associated generalised series $\tilde{y}_0=y_0d_k^n$. 

 If $k<k_0$, by the Corollary \ref{coro:k_0}, then  $\tilde{k}>k$. It implies that \begin{center}
$\begin{array}{lcl}
\tilde{y}_0&=&y_{0,k}d_k^n+\cdots+y_{0,1}d_k^n\\
&=&\tilde{y}_{0,k}+\cdots+\tilde{y}_{0,1}\in\mathds{K}_r^\prec
\end{array}$ 
\end{center} as required to apply the Lemma \ref{lemme:princ}.  Note that the support of $y_0$ is then equal to $\textrm{Supp }\tilde{y}_0 -nv(d_k)$. Therefore, it is included into $\textrm{Supp }\tilde{y}_0+ \mathbb{N}(-v(d_k))$, which is the application of an elementary transformation to $\textrm{Supp }\tilde{y}_0$.

If $k=k_0$, since
 $v(y_0^{(i)})>\underline{0}$ for any $i=0,\ldots,n$, then $v(y_0)$ verifies the Condition \ref{eq:cond-posit}. So, $v(\tilde{y}_0)$ is also positive, and $\tilde{y}_0$ is of the form \begin{center}
$\tilde{y}_0=\tilde{y}_{0,l}+\cdots+\tilde{y}_{0,1}\in\mathds{K}_r^\prec$
\end{center}
with $l\geq k_0$. Then, if necessary, we apply in the differential series $\tilde{F}_{0,k_0}$ just obtained, the change of derivation from $D_{k_0}$ to $D_l$. Then we can apply to the latter, together with the series $\tilde{y}_0$, the Lemma \ref{lemme:princ}. As before, $\textrm{Supp }y_0$ is deduced from $\textrm{Supp }\tilde{y}_0$ by an elementary transformation.\\

In any case, we are reduced to the hypothesis of the Lemma \ref{lemme:princ}. To simplify the notations, let us denote in any of the described cases, the obtained differential series with Weierstrass order $w$ by $G(y,\ldots,D_k^ny)$, and the corresponding  generalised series by $y_0=y_{0,k}+\cdots+y_{0,1}$. Applying the Lemma \ref{lemme:princ}, we obtain a dichotomy for $y_{0,k}$, the second case being exactly the second one in the statement of Theorem \ref{theo_princ}.

For the first one, we have $\textrm{Supp}\ y_{0,k}\subset\mathcal{R}_k$, where $\mathcal{R}_k$ is a well-ordered subset of $\Gamma_{>0}$ obtained from $\textrm{Supp}\ G$ (and therefore from $\textrm{Supp}\ F$) and  $\mathcal{T}_k,\ldots,\mathcal{T}_r$ by a finite number of elementary transformations. Then we perform the additive conjugation $y=\tilde{y}+y_{0,k}$ in the original differential series $F$, and we get a new one $F_{\textrm{Supp}\ y_{0,k}}(\tilde{y},\ldots,\tilde{y}^{(n)})$ together with a series $\tilde{y}_0=y_{0,l}+\ldots+y_{0,1}$ for some $l<k$. Moreover, by the Proposition \ref{transl:propos}, the support of $F_{\textrm{Supp}\ y_{0,k}}$ is obtained from the one of $F$ and the sets $\mathcal{T}_k$ and $\textrm{Supp}\ y_{0,k}$ by finitely many elementary transformations. But $\textrm{Supp}\ y_{0,k}\subset\mathcal{R}_k$ is obtained itself from $\textrm{Supp}\ F,\ \mathcal{T}_k,\ldots,\mathcal{T}_r$ by a finite number of elementary transformations. Then, so $\textrm{Supp}\ F_{\textrm{Supp}\ y_{0,k}}$ is. We can resume the preceding arguments with $F_{\textrm{Supp}\ y_{0,k}}(\tilde{y},\ldots,\tilde{y}^{(n)})$ and $\tilde{y}_0$.

Thus we prove gradually in the case where there is no stabilisation that $\textrm{Supp}\ y_{0}\subset\mathcal{R}=\sum_{l=1}^k\left\langle\mathcal{R}_l\right\rangle$, the latter being obtained from $\textrm{Supp}\ F$, $\mathcal{T}_k,\ldots,\mathcal{T}_r$ by finitely many elementary transformations as desired.\hfill$\Box$
\begin{center}

\end{center}
Concerning the proof of our main lemma, we consider the set $\{1,\ldots,r\}$ endowed with the \textsl{reverse ordering} $1>2>\cdots>r$, and we denote it $\{r,\ldots,1\}$. Then we define the lexicographical product $\overrightarrow{\mathbb{N}\times\{r,\ldots,1\}}$, remark that it is a well-ordered set, and denote $(w,k)$ its elements. To prove Lemma \ref{lemme:princ}, we proceed by transfinite induction on $(w,k)\in\overrightarrow{\mathbb{N}\times\{r,\ldots,1\}}$.

First, if the Differential Series \ref{eq:serie-diff-D_k} has Weierstrass order 0, then for any $k\in\{r,\ldots,1\}$ and any initial part $p$ of $y_0$, \begin{center}
$F(p,\ldots,D_k^np)=c_0+\sum_{I\in\mathbb{N}^{n+1}\backslash \{0\}}
c_Ip^{(I)}$.
\end{center} But for any $I\in\mathbb{N}^{n+1}$ with $|I|\geq 1$, $v(y_0^{(I)_k})>\underline{0}$ (Lemma \ref{deriv_posit}) and by hypothesis $v(c_0)=\underline{0}\leq v(c_I)$. So 
$v(F(p,\ldots,D_k^np)) =v(c_0)=\underline{0}$, which means that $y_{0}$ stabilises on $F$ with initial part 0.

Second, we consider $(w,k)\in\overrightarrow{\mathbb{N}^*\times\{r,\ldots,1\}}$ and suppose that the Theorem \ref{theo_princ} holds for any $(\tilde{w},\tilde{k})<(w,k)$. We treat separately the case $w=1$.

\subsection{The case $w=1$.}\label{sect:weier1}
\begin{notation}Given a Differential Series \ref{eq:serie-diff-D_k} with Weierstrass order $1$, we denote:\begin{itemize}
\item $\mathcal{A}$ the set of all multi-index $I$ with length 1 such that $v(c_I)=\underline{0}$;
\item for all $I\in\mathbb{N}^{n+1}$, $c_{I,0}$ the leading coefficient of $c_I$.
\end{itemize}
\end{notation}

\begin{defi}\label{defi:poly_annul}
With the terminology of \cite{ince:ODE}, we call \textsl{indicial polynomial} associated to a Weierstrass order 1 series $F(y,\ldots,D_k^ny)$ the polynomial $\pi(X)=\sum_{I\in\mathcal{A}}c_{I,0}X^{\|I\|}$.\\
We denote $\{\rho_1,\ldots,\rho_m\}$ the set of its real positive roots ($m\leq n$ since its degree is at most $n$).
\end{defi}

\begin{lemma}\label{lemma:dichot_weier1}
Given a Weierstrass order 1 Differential Series \ref{eq:serie-diff-D_k} and a series $y_0\in\mathds{K}_r^\prec$ with initial part \begin{center}
$y_{0,k}=\sum_{\mu\in\textrm{Supp}\ y_{0,k}}m_\mu t^\mu\in\mathds{K}_{r,k}$,
\end{center} there are two cases:\\
\begin{enumerate}
\item either for any proper initial segment $S$ of $\textrm{Supp}\ y_{0,k}$ with successor element $\mu=(0,\ldots,0,\mu_k,\ldots,\mu_r)$ ($\mu_k>0$), we have $\mu=v_S$, or $\mu_k\in\{\rho_1,\ldots,\rho_m\}$;
\item or there exists a proper initial segment $S_0$ of $\textrm{Supp}\ y_{0,k}$ with successor element $\mu^{(0)}$ such that the two following properties hold:
\begin{enumerate}
\item for any successor segment $\tilde{S}$ of $S_0$ in $\textrm{Supp}\ y_0$, we have \begin{center}
$v_{S_0\cup \tilde{S}}=v_{S_0\cup \{\mu^{(0)}\}}=\min\{v_{S_0},\mu^{(0)}\}\leq\mu^{(0)}$;
\end{center}
\item for any proper initial segment $S\subsetneq S_0$ with successor element $\mu<\mu^{(0)}$, we have $\mu=v_S$ or $\mu_k\in\{\rho_1,\ldots,\rho_m\}$.
\end{enumerate}
\end{enumerate}
\end{lemma}
\begin{demo}
We consider a proper initial segment $S$ of $\textrm{Supp}\ y_{0,k}$ with successor element $\mu$, and some successor segment $\tilde{S}$ of $S$ (so $\mu$ is the least element of $\tilde{S}$). From the Taylor expansion formula \ref{taylor}: \begin{center}
$f_{S\cup\tilde{S}}=\sum_{I\in\mathbb{N}^{n+1}} (f_S^{(I)}/I!)p_{\tilde{S}}^{(I)}$.
\end{center}
But $p_{\tilde{S}}=m_\mu t^\mu(1+\epsilon)$ and so \begin{center}
$p_{\tilde{S}}^{(I)}=m_\mu^{|I|}\mu_k^{\|I\|}t^{|I|\mu}(1+\epsilon_I)$ 
\end{center}for some $\epsilon,\epsilon_I\in\mathds{K}_r^\prec$ (Proposition \ref{propo_D_k}). Moreover, for any $I\in\mathcal{A}$, $\delta((f_S^{(I)}/I!))=c_{I,0}$. So \begin{center}
$\delta((f_S^{(I)}/I!)p_{\tilde{S}}^{(I)})=m_\mu t^{\mu}\mu_k^{\|I\|}c_{I,0}$.
\end{center} For any $J\in\mathbb{N}^{n+1}\backslash\{0\}$, $J\notin\mathcal{A}$, we have \begin{center}
$v((f_S^{(J)}/J!)p_{\tilde{S}}^{(J)})=v_S^{(J)}+|J|\mu>\mu$. 
\end{center} Thus we obtain $f_{S\cup\tilde{S}}=f_S+m_\mu t^\mu(\pi(\mu_k)+\epsilon)$ for some $\epsilon\in\mathds{K}_r^\prec$. Then the dichotomy of the lemma follows from the ultrametric triangular inequality for the valuation $v$. For the second case, $S_0$ is the least initial segment such that $\mu^{(0)}\neq v_{S_0}$ and $\mu_k^{(0)}\notin\{\rho_1,\ldots,\rho_m\}$. Then for any successor segment $\tilde{S}$ of $S_0$, we have \begin{center}
$v_{S_0\cup \tilde{S}}=v_{S_0\cup \{\mu^{(0)}\}}=\min\{v_{S_0},\mu^{(0)}\}$.
\end{center}
\end{demo}

Returning to the proof of Lemma \ref{lemme:princ}, in case $k=r$, we set \begin{center}
$\mathcal{R}=\langle\textrm{Supp}\ F\rangle+ \langle(0,\ldots,0,\rho_1),\ldots,(0,\ldots,0,\rho_m)\rangle+ \mathcal{T}_r$.
\end{center} Then, from Proposition \ref{supp_eval}, we remark that \begin{center}
$v_S\in \textrm{Supp}\
F+\mathcal{T}_r+\langle S\rangle$
\end{center} for any initial segment $S$ of $\textrm{Supp}\ y_{0,r}$. So using the relations $\mu=v_S$ or $\mu_r\in\{\rho_1,\ldots,\rho_m\}$, by a straightforward transfinite induction, we obtain that $\textrm{Supp}\ y_0\subset\mathcal{R}$ in case (1) of the preceding lemma, respectively $S_0\subset\mathcal{R}$ in case (2). The subcase (2)(a) means exactly that $y_{0}$ stabilises on $F$ with initial part $p_{S_0\cup \{\mu^{(0)}\}}$.

In case $k\in\{r-1,\ldots,1\}$, we consider $y_0=y_{0,k}+\cdots+y_{0,1}$ and a Differential Series \ref{eq:serie-diff-D_k} with Weierstrass order 1. We suppose that the Lemma \ref{lemme:princ} holds for any $l\in\{r,\ldots,k+1\}$.
According to the dichotomy in the Lemma \ref{lemma:dichot_weier1}, we show by transfinite induction that $\textrm{Supp}\ y_{0,k}$ (respectively $S_0$) is included in an additive sub-semigroup $\mathcal{R}_k$ of $\Gamma_{>0}$ of type \begin{center}
$\mathcal{R}_k=\left\langle \textrm{Supp}\ F\right\rangle+\mathcal{T}_k+ \cdots+\mathcal{T}_r+\mathcal{R}_{1,k}+\cdots+\mathcal{R}_{m,k}$
\end{center} where $m$ is the number of positive roots of the indicial polynomial associated to $F$ (see Definition \ref{defi:poly_annul}), and the $\mathcal{R}_{i,k}$'s are additive sub-semigroups of $\Gamma_{>0}$ obtained from $\textrm{Supp}\ F$, $\mathcal{T}_k,\ldots,\mathcal{T}_r$ by a finite number of elementary transformations.

Indeed, let us consider a proper initial segment $S$ of $\textrm{Supp}\ y_{0,k}$ (respectively $S_0$) with successor element $\mu=(0,\ldots,0,\mu_k,\ldots,\mu_r)$. We suppose that $S$ is included in an additive semi-group $\mathcal{R}$ that contains $\left\langle \textrm{Supp}\ F\right\rangle,\ \mathcal{T}_k,\ldots,\mathcal{T}_r$. From Lemma \ref{lemma:dichot_weier1}, there are two cases. Either $\mu=v_S$. But from Proposition \ref{supp_eval} we have \begin{center}
$v_S\in\textrm{Supp}\ f_S\subset\textrm{Supp}\ F+\mathcal{T}_k+\left\langle S\right\rangle$
\end{center} and $\left\langle S\right\rangle\subset\mathcal{R}$ as well as $\textrm{Supp}\ F$ and $\mathcal{T}_k$. So $\mu\in\mathcal{R}$.

 Or $\mu_k\in\{\rho_1,\ldots,\rho_m\}$. For instance, $\mu_k=\rho_h$ for some $h\in\{1,\ldots,m\}$ fixed.
We set $S_h=\{\mu\in\textrm{Supp}\ y_{0,k}\ |\ \mu_k=\rho_h\}$, $\sigma^{(h)}=\min S_h$, \begin{center}
$S'_h=\{\mu\in\textrm{Supp}\ y_{0,k}\ |\ \mu<\sigma^{(h)}\}$
\end{center} and $S''_h=S'_h\cup\{\sigma^{(h)}\}$. So $S'_h\subset S\subset\mathcal{R}$ and $S''_h\subset\mathcal{R}+\mathbb{N}\sigma^{(h)}$. Let us show that there exists an additive sub-semigroup $\mathcal{R}_{h,k}$ of $\Gamma_{>0}$ obtained from $\mathcal{R}$, $\textrm{Supp}\ F$, $\mathcal{T}_k,\ldots,\mathcal{T}_r$ by a finite number of elementary transformations such that $S_h\subset \mathcal{R}_{h,k}$.

We notice that \begin{center}
$\sigma^{(h)}=(0,\ldots,0,\rho_h,0,\ldots,0,\sigma_{l,h},\ldots,\sigma_{r,h})$ 
\end{center}for some $l\in\{k+1,\ldots,r\}$ with $\sigma_{l,h}\neq 0$. So  \begin{center}
$p_{S_h\backslash\{\sigma^{(h)}\}}/t^{\sigma^{(h)}}= z_{0,l}+\cdots+z_{0,k+1}\in\mathds{K}_r^\prec$
\end{center} with initial part $z_{0,l}\in\mathds{K}_{r,l}$. Then we reduce to a differential series $G(z,\ldots,D_l^nz)$ together with the series $z_{0,l}\in\mathds{K}_{r,l}$ so as to apply the induction hypothesis.

The calculations we make are represented symbolically in cases $r=2$ and $k=1$ by the following picture, in which the black points represent elements of the support of $y_{0,k}$. The necessary changes of derivation are also mentioned.\\

\begin{figure}
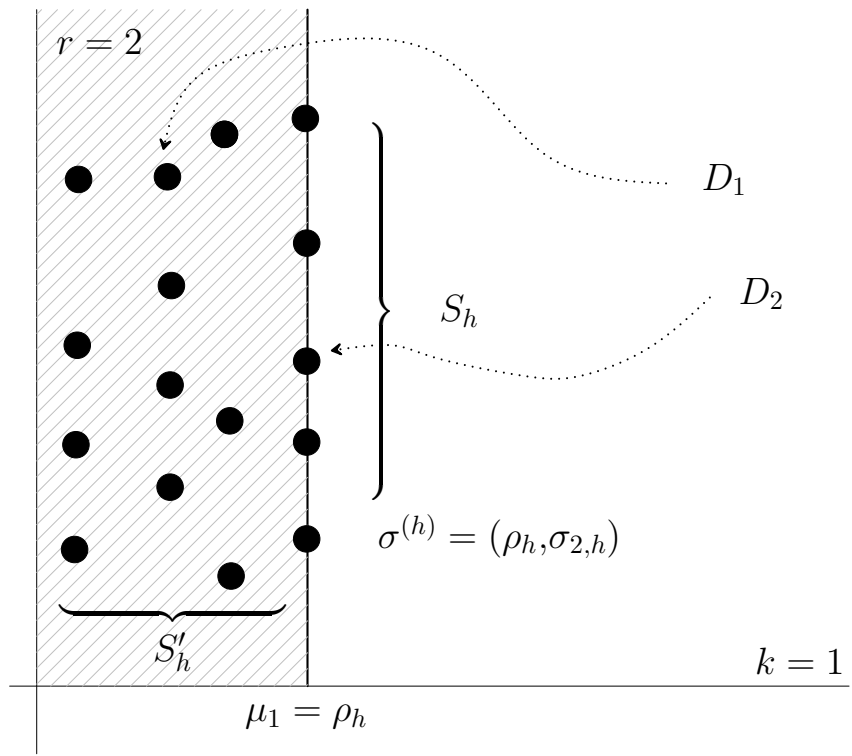
\label{fig:support_y_0}
\begin{pgfpicture}{-0.5cm}{-0.5cm}{11.66cm}{10.4cm}%
\pgfsetroundjoin \pgfsetroundcap%
\pgfsetlinewidth{0.2pt} 
\pgfxyline(0,0.9)(11.16,0.9)\pgfxyline(0.36,0)(0.36,9.9)
\pgfsetlinewidth{0.8pt} 
\pgfxyline(3.96,0.9)(3.96,9.9)
\pgfsetlinewidth{0.2pt}\pgfsetstrokecolor{rgb,1:red,0.7529;green,0.7529;blue,0.7529}
\pgfxyline(0.4243,9.9)(0.36,9.8357)\pgfxyline(0.6364,9.9)(0.36,9.6236)\pgfxyline(0.8485,9.9)(0.36,9.4115)
\pgfxyline(1.0607,9.9)(0.36,9.1993)\pgfxyline(1.2728,9.9)(0.36,8.9872)\pgfxyline(1.4849,9.9)(0.36,8.7751)
\pgfxyline(1.6971,9.9)(0.36,8.5629)\pgfxyline(1.9092,9.9)(0.36,8.3508)\pgfxyline(2.1213,9.9)(0.36,8.1387)
\pgfxyline(2.3335,9.9)(0.36,7.9265)\pgfxyline(2.5456,9.9)(0.36,7.7144)\pgfxyline(2.7577,9.9)(0.36,7.5023)
\pgfxyline(2.9698,9.9)(0.36,7.2902)\pgfxyline(3.182,9.9)(0.36,7.078)\pgfxyline(3.3941,9.9)(0.36,6.8659)
\pgfxyline(3.6062,9.9)(0.36,6.6538)\pgfxyline(3.8184,9.9)(0.36,6.4416)\pgfxyline(3.96,9.8295)(0.36,6.2295)
\pgfxyline(3.96,9.6174)(0.36,6.0174)\pgfxyline(3.96,9.4052)(0.36,5.8052)\pgfxyline(3.96,9.1931)(0.36,5.5931)
\pgfxyline(3.96,8.981)(0.36,5.381)\pgfxyline(3.96,8.7688)(0.36,5.1688)\pgfxyline(3.96,8.5567)(0.36,4.9567)
\pgfxyline(3.96,8.3446)(0.36,4.7446)\pgfxyline(3.96,8.1324)(0.36,4.5324)\pgfxyline(3.96,7.9203)(0.36,4.3203)
\pgfxyline(3.96,7.7082)(0.36,4.1082)\pgfxyline(3.96,7.496)(0.36,3.896)\pgfxyline(3.96,7.2839)(0.36,3.6839)
\pgfxyline(3.96,7.0718)(0.36,3.4718)\pgfxyline(3.96,6.8596)(0.36,3.2596)\pgfxyline(3.96,6.6475)(0.36,3.0475)
\pgfxyline(3.96,6.4354)(0.36,2.8354)\pgfxyline(3.96,6.2232)(0.36,2.6232)\pgfxyline(3.96,6.0111)(0.36,2.4111)
\pgfxyline(3.96,5.799)(0.36,2.199)\pgfxyline(3.96,5.5869)(0.36,1.9869)\pgfxyline(3.96,5.3747)(0.36,1.7747)
\pgfxyline(3.96,5.1626)(0.36,1.5626)\pgfxyline(3.96,4.9505)(0.36,1.3505)\pgfxyline(3.96,4.7383)(0.36,1.1383)
\pgfxyline(3.96,4.5262)(0.36,0.9262)\pgfxyline(3.96,4.3141)(0.5459,0.9)\pgfxyline(3.96,4.1019)(0.7581,0.9)
\pgfxyline(3.96,3.8898)(0.9702,0.9)\pgfxyline(3.96,3.6777)(1.1823,0.9)\pgfxyline(3.96,3.4655)(1.3945,0.9)
\pgfxyline(3.96,3.2534)(1.6066,0.9)\pgfxyline(3.96,3.0413)(1.8187,0.9)\pgfxyline(3.96,2.8291)(2.0309,0.9)
\pgfxyline(3.96,2.617)(2.243,0.9)\pgfxyline(3.96,2.4049)(2.4551,0.9)\pgfxyline(3.96,2.1927)(2.6673,0.9)
\pgfxyline(3.96,1.9806)(2.8794,0.9)\pgfxyline(3.96,1.7685)(3.0915,0.9)\pgfxyline(3.96,1.5563)(3.3037,0.9)
\pgfxyline(3.96,1.3442)(3.5158,0.9)\pgfxyline(3.96,1.1321)(3.7279,0.9)\pgfxyline(3.96,0.9199)(3.9401,0.9)
\pgfsetstrokecolor{black}
\pgfputat{\pgfxy(3.9492,0.4941)}{\pgfnode{rectangle}{center}{\color{black}\Large $\mu_1=\rho_{h}$}{}{\pgfusepath{}}}
\pgfputat{\pgfxy(6.5938,2.894)}{\pgfnode{rectangle}{center}{\color{black}\Large $\sigma^{\left(h\right)}=\left(\rho_h,\sigma_{2,h}\right)$ }{}{\pgfusepath{}}}
\pgfputat{\pgfxy(10.5077,1.2353)}{\pgfnode{rectangle}{center}{\color{black}\Large $k=1$}{}{\pgfusepath{}}}
\pgfputat{\pgfxy(1.1988,9.4765)}{\pgfnode{rectangle}{center}{\color{black}\Large $r=2$}{}{\pgfusepath{}}}
\pgfsetfillcolor{black}\pgfellipse[fillstroke]{\pgfxy(3.9492,2.8588)}{\pgfxy(0.1762,0)}{\pgfxy(0,0.1762)}
\pgfellipse[fillstroke]{\pgfxy(3.9492,4.147)}{\pgfxy(0.1762,0)}{\pgfxy(0,0.1762)}
\pgfellipse[fillstroke]{\pgfxy(3.9492,5.2236)}{\pgfxy(0.1762,0)}{\pgfxy(0,0.1762)}
\pgfellipse[fillstroke]{\pgfxy(3.9492,6.7941)}{\pgfxy(0.1762,0)}{\pgfxy(0,0.1762)}
\pgfellipse[fillstroke]{\pgfxy(3.9316,8.453)}{\pgfxy(0.1762,0)}{\pgfxy(0,0.1762)}
\pgfputat{\pgfxy(5.9767,5.8941)}{\pgfnode{rectangle}{center}{\Large $S_h$}{}{\pgfusepath{}}}
{\pgftransformrotate{-90}
\pgfputat{\pgfxy(4.9365,5.8941)}{\pgfnode{rectangle}{center}{\large $\overbrace{\hspace{5cm}}$}{}{\pgfusepath{}}}}
\pgfellipse[fillstroke]{\pgfxy(0.8638,2.7176)}{\pgfxy(0.1762,0)}{\pgfxy(0,0.1762)}
\pgfellipse[fillstroke]{\pgfxy(0.8815,4.1117)}{\pgfxy(0.1762,0)}{\pgfxy(0,0.1762)}
\pgfellipse[fillstroke]{\pgfxy(0.8991,5.4353)}{\pgfxy(0.1762,0)}{\pgfxy(0,0.1762)}
\pgfellipse[fillstroke]{\pgfxy(0.9167,7.6412)}{\pgfxy(0.1762,0)}{\pgfxy(0,0.1762)}
\pgfellipse[fillstroke]{\pgfxy(2.0981,7.6765)}{\pgfxy(0.1762,0)}{\pgfxy(0,0.1762)}
\pgfellipse[fillstroke]{\pgfxy(2.1508,6.2294)}{\pgfxy(0.1762,0)}{\pgfxy(0,0.1762)}
\pgfellipse[fillstroke]{\pgfxy(2.1334,4.9059)}{\pgfxy(0.1762,0)}{\pgfxy(0,0.1762)}
\pgfellipse[fillstroke]{\pgfxy(2.1334,3.5471)}{\pgfxy(0.1762,0)}{\pgfxy(0,0.1762)}
\pgfellipse[fillstroke]{\pgfxy(2.9443,2.3647)}{\pgfxy(0.1762,0)}{\pgfxy(0,0.1762)}
\pgfellipse[fillstroke]{\pgfxy(2.9266,4.4294)}{\pgfxy(0.1762,0)}{\pgfxy(0,0.1762)}
\pgfellipse[fillstroke]{\pgfxy(2.8561,8.2411)}{\pgfxy(0.1762,0)}{\pgfxy(0,0.1762)}
\pgfellipse[fillstroke]{\pgfxy(2.8561,8.2411)}{\pgfxy(0.1762,0)}{\pgfxy(0,0.1762)}
{\pgftransformrotate{-180}
\pgfputat{\pgfxy(2.1685,1.8706)}{\pgfnode{rectangle}{center}{\large $\overbrace{\hspace{3cm}}$}{}{\pgfusepath{}}}}
\pgfputat{\pgfxy(2.1685,1.3588)}{\pgfnode{rectangle}{center}{\Large $S'_h$}{}{\pgfusepath{}}}
\pgfsetlinewidth{0.8pt} 
\pgfsetdash{{0pt}{3pt}}{0pt}
\pgfmoveto{\pgfxy(8.7269,7.5883)}\pgflineto{\pgfxy(8.5685,7.5958)}\pgflineto{\pgfxy(8.4105,7.6049)}\pgflineto{\pgfxy(8.2535,7.6169)}\pgflineto{\pgfxy(8.0979,7.6333)}\pgflineto{\pgfxy(7.9441,7.6555)}
\pgflineto{\pgfxy(7.7927,7.6851)}\pgflineto{\pgfxy(7.6441,7.7234)}\pgflineto{\pgfxy(7.4988,7.7721)}\pgflineto{\pgfxy(7.3573,7.8324)}\pgflineto{\pgfxy(7.2199,7.906)}\pgflineto{\pgfxy(7.0873,7.9942)}
\pgflineto{\pgfxy(6.9596,8.0977)}\pgflineto{\pgfxy(6.8356,8.2142)}\pgflineto{\pgfxy(6.7143,8.3404)}\pgflineto{\pgfxy(6.5941,8.4732)}\pgflineto{\pgfxy(6.4739,8.6093)}\pgflineto{\pgfxy(6.3522,8.7456)}
\pgflineto{\pgfxy(6.2279,8.8788)}\pgflineto{\pgfxy(6.0995,9.0059)}\pgflineto{\pgfxy(5.9658,9.1235)}\pgflineto{\pgfxy(5.8254,9.2285)}\pgflineto{\pgfxy(5.677,9.3177)}\pgflineto{\pgfxy(5.5199,9.3887)}
\pgflineto{\pgfxy(5.3549,9.4424)}\pgflineto{\pgfxy(5.1838,9.4805)}\pgflineto{\pgfxy(5.0079,9.5046)}\pgflineto{\pgfxy(4.8288,9.5165)}\pgflineto{\pgfxy(4.6481,9.5177)}\pgflineto{\pgfxy(4.4673,9.51)}
\pgflineto{\pgfxy(4.2878,9.4951)}\pgflineto{\pgfxy(4.1113,9.4746)}\pgflineto{\pgfxy(3.9391,9.4502)}\pgflineto{\pgfxy(3.773,9.4235)}\pgflineto{\pgfxy(3.6141,9.396)}\pgflineto{\pgfxy(3.4628,9.3675)}
\pgflineto{\pgfxy(3.3193,9.3374)}\pgflineto{\pgfxy(3.1837,9.3053)}\pgflineto{\pgfxy(3.0562,9.2707)}\pgflineto{\pgfxy(2.9368,9.2331)}\pgflineto{\pgfxy(2.8258,9.1919)}\pgflineto{\pgfxy(2.7233,9.1468)}
\pgflineto{\pgfxy(2.6294,9.0972)}\pgflineto{\pgfxy(2.5443,9.0425)}\pgflineto{\pgfxy(2.4682,8.9824)}\pgflineto{\pgfxy(2.4009,8.9164)}\pgflineto{\pgfxy(2.3418,8.845)}\pgflineto{\pgfxy(2.2902,8.7687)}
\pgflineto{\pgfxy(2.2452,8.688)}\pgflineto{\pgfxy(2.2059,8.6035)}\pgflineto{\pgfxy(2.1715,8.5158)}\pgflineto{\pgfxy(2.1413,8.4253)}\pgflineto{\pgfxy(2.1144,8.3327)}\pgflineto{\pgfxy(2.0899,8.2384)}
\pgflineto{\pgfxy(2.0671,8.143)}\pgflineto{\pgfxy(2.0452,8.0471)}\pgfstroke
\pgfsetdash{}{0pt}
\pgfmoveto{\pgfxy(2.1359,8.1447)}\pgflineto{\pgfxy(2.0452,8.0471)}\pgflineto{\pgfxy(2.006,8.1745)}\pgflineto{\pgfxy(2.0581,8.1034)}\pgfclosepath\pgffillstroke
\pgfsetdash{{0pt}{3pt}}{0pt}
\pgfmoveto{\pgfxy(9.3089,6.0705)}\pgflineto{\pgfxy(9.2389,6.0043)}\pgflineto{\pgfxy(9.1688,5.9384)}\pgflineto{\pgfxy(9.0983,5.873)}\pgflineto{\pgfxy(9.0272,5.8083)}\pgflineto{\pgfxy(8.9554,5.7448)}
\pgflineto{\pgfxy(8.8826,5.6825)}\pgflineto{\pgfxy(8.8087,5.6218)}\pgflineto{\pgfxy(8.7334,5.563)}\pgflineto{\pgfxy(8.6566,5.5062)}\pgflineto{\pgfxy(8.5782,5.4518)}\pgflineto{\pgfxy(8.4978,5.4)}
\pgflineto{\pgfxy(8.4155,5.3511)}\pgflineto{\pgfxy(8.3316,5.305)}\pgflineto{\pgfxy(8.2467,5.262)}\pgflineto{\pgfxy(8.1615,5.2219)}\pgflineto{\pgfxy(8.0763,5.185)}\pgflineto{\pgfxy(7.9918,5.151)}
\pgflineto{\pgfxy(7.9085,5.1203)}\pgflineto{\pgfxy(7.827,5.0927)}\pgflineto{\pgfxy(7.7479,5.0683)}\pgflineto{\pgfxy(7.6716,5.0472)}\pgflineto{\pgfxy(7.5987,5.0294)}\pgflineto{\pgfxy(7.5295,5.0149)}
\pgflineto{\pgfxy(7.4633,5.0036)}\pgflineto{\pgfxy(7.3989,4.9955)}\pgflineto{\pgfxy(7.3353,4.9903)}\pgflineto{\pgfxy(7.2714,4.9881)}\pgflineto{\pgfxy(7.2062,4.9887)}\pgflineto{\pgfxy(7.1386,4.9919)}
\pgflineto{\pgfxy(7.0675,4.9977)}\pgflineto{\pgfxy(6.9919,5.006)}\pgflineto{\pgfxy(6.9107,5.0165)}\pgflineto{\pgfxy(6.8229,5.0294)}\pgflineto{\pgfxy(6.7278,5.0443)}\pgflineto{\pgfxy(6.6265,5.0611)}
\pgflineto{\pgfxy(6.5206,5.0796)}\pgflineto{\pgfxy(6.4117,5.0995)}\pgflineto{\pgfxy(6.3013,5.1206)}\pgflineto{\pgfxy(6.191,5.1426)}\pgflineto{\pgfxy(6.0822,5.1654)}\pgflineto{\pgfxy(5.9767,5.1886)}
\pgflineto{\pgfxy(5.8759,5.2121)}\pgflineto{\pgfxy(5.7813,5.2356)}\pgflineto{\pgfxy(5.6947,5.2589)}\pgflineto{\pgfxy(5.6169,5.2817)}\pgflineto{\pgfxy(5.5473,5.3039)}\pgflineto{\pgfxy(5.4845,5.3251)}
\pgflineto{\pgfxy(5.4274,5.3452)}\pgflineto{\pgfxy(5.3745,5.364)}\pgflineto{\pgfxy(5.3248,5.3812)}\pgflineto{\pgfxy(5.2768,5.3965)}\pgflineto{\pgfxy(5.2294,5.4098)}\pgflineto{\pgfxy(5.1812,5.4209)}
\pgflineto{\pgfxy(5.1311,5.4294)}\pgflineto{\pgfxy(5.0776,5.4353)}\pgflineto{\pgfxy(5.0199,5.4382)}\pgflineto{\pgfxy(4.9581,5.4385)}\pgflineto{\pgfxy(4.8926,5.4364)}\pgflineto{\pgfxy(4.8237,5.4321)}
\pgflineto{\pgfxy(4.752,5.426)}\pgflineto{\pgfxy(4.6778,5.4182)}\pgflineto{\pgfxy(4.6016,5.4091)}\pgflineto{\pgfxy(4.5237,5.399)}\pgflineto{\pgfxy(4.4445,5.388)}\pgflineto{\pgfxy(4.3646,5.3765)}
\pgflineto{\pgfxy(4.2842,5.3647)}\pgfstroke
\pgfsetdash{}{0pt}
\pgfmoveto{\pgfxy(4.4081,5.3155)}\pgflineto{\pgfxy(4.2842,5.3647)}\pgflineto{\pgfxy(4.3888,5.4474)}\pgflineto{\pgfxy(4.3413,5.3731)}\pgfclosepath\pgffillstroke
\pgfputat{\pgfxy(9.5027,7.6412)}{\pgfnode{rectangle}{center}{\Large $D_1$}{}{\pgfusepath{}}}
\pgfputat{\pgfxy(9.9965,6.1236)}{\pgfnode{rectangle}{center}{\Large $D_2$}{}{\pgfusepath{}}}
\end{pgfpicture}
\caption{Representation of the support of $y_0$, with $r=2$ and $k=1$}
\end{figure}

\noindent(1) We perform the \textsl{additive conjugation} $y=\tilde{y}+p_{S''_h}$ in the differential series\\ $F(y,\ldots,D_k^ny)$. From Proposition \ref{transl:propos} and Remark \ref{rem:weier_transl}, we obtain a new series \begin{center}
$F_{S''_h}(\tilde{y},\ldots,D_k^n\tilde{y})= \sum_{I\in\mathbb{N}^{n+1}}c'_I\tilde{y}^{(I)}$ 
\end{center}with Weierstrass order 1. We have $c'_I=f_{S''_h}^{(I)}/I!$ and in particular, for any $I\in\mathcal{A}$, $v(c'_I)=v(c_I)=\underline{0}$. Moreover \begin{center}
$\textrm{Supp}\ F_{S''_h}\subset\textrm{Supp}\ F+\mathcal{T}_k+\left\langle S''_h\right\rangle\subset\mathcal{R}+\mathbb{N}\sigma^{(h)}$.
\end{center} The associated series is $\tilde{y}_0=y_0-p_{S''_h}$ with initial part $\tilde{y}_{0,k}=y_{0,k}-p_{S''_h}\in\mathds{K}_{r,k}$.\\

\noindent(2) We perform in $F_{S''_h}$ the \textsl{multiplicative conjugation} $y_1=t^{\sigma^{(h)}}z$. From Proposition \ref{eclt:propos}, we obtain a differential series \begin{center}
$\hat{F}(z,\ldots,D_k^nz)=\sum_{J\in\mathbb{N}^{n+1}}\hat{c}_Jz^{(J)}$ 
\end{center}with \begin{center}
$v(\hat{c}_J)\geq\min\{v(c'_I)+v((t^{\sigma^{(h)}})^{(K)})\ \mid\ |I|=|J|=|K|,\ \|I\|=\|K\|+\|J\|,\ J,K\leq_{antilex}I\}$
\end{center} and \begin{center}
$(t^{\sigma^{(h)}})^{(K)}\sim\rho_h^{\|K\|}t^{|K|\sigma^{(h)}}= \rho_h^{\|K\|}t^{|I|\sigma^{(h)}}$.
\end{center} So \begin{center}
$v(\hat{c}_J)\geq\min\{v(c'_I)+ |I|\sigma^{(h)}\ \mid\ |I|=|J|,\ J\leq_{antilex}I\}$.
\end{center} According to the Remark \ref{rem:weier_eclt}, we have $v(\hat{c}_{I_0})=\sigma^{(h)}$ where $I_0$ is the greatest multi-index from $\mathcal{A}$ for anti-lexicographical ordering (in particular $|I_0|=1$). Moreover, from Proposition \ref{eclt:propos}, \begin{center}
$\textrm{Supp}\ \hat{F}\subset\textrm{Supp}\ F_{S''_h}+\mathcal{T}_k+\mathbb{N}\sigma^{(h)}\subset \mathcal{R}+\mathbb{N}\sigma^{(h)}$.
\end{center} 
The associated series is $z_0=\tilde{y}_0/t^{\sigma^{(h)}}$, which has $z_{0,l}\in\mathds{K}_{r,l}$  as initial part.\\

 \noindent(3) We perform the \textsl{change of derivation} \begin{center}
$\hat{F}(z,\ldots,D_k^nz)=\tilde{F}_{k,l}(z,\ldots,D_l^nz)= \sum_{L\in\mathbb{N}^{n+1}}\tilde{c}_Lz^{(L)}$
\end{center} as in the Proposition \ref{chgt_cond1}. From the Remark \ref{rem:chge_deriv}, the terms $\hat{c}_Jz^{(J)}$ with $v(\hat{c}_J)=\sigma^{(h)}$ provide at least one term $\tilde{c}_{L}z^{(L)}$ with $|L|=1$ and \begin{center}
$v(\tilde{c}_{L})=(0,\ldots,0,\rho_h,\chi_{k+1},\ldots,\chi_r)$
\end{center} for some reals $\chi_{k+1},\ldots,\chi_r$ ($L=J_0$ being the greatest element for anti-lexicographical ordering among these $J$'s is such a good candidate). Then we observe that \begin{center}
$\min\{v(\tilde{c}_L), L\in\mathbb{N}^{n+1}, |L|\geq 1\}=v(\tilde{c}_{L_0})=(0,\ldots,0,\rho_h,\chi^{(0)}_{k+1},\ldots,\chi^{(0)}_r)$ 
\end{center}
 for some $L_0$ with $|L_0|=1$. We set $v^{(0)}=v(\tilde{c}_{L_0})$ and there are two cases. Either $v(\tilde{c}_0)\leq v^{(0)}$. Then for any initial segment $\tilde{S}$ of $\textrm{Supp}\ z_{0,l}$, we have $v(\tilde{f}_{\tilde{S}})=v(\tilde{c}_0)$. It means that $z_0$ stabilises on $\tilde{F}_{k,l}$ with initial part 0, and equivalently $y_0$ stabilises on $F$ with initial part $p_{S''_h}$.\\ 
Or $v(\tilde{c}_0)> v^{(0)}$. Then we divide $\tilde{F}_{k,l}(z,\ldots,D_l^nz)$ by $t^{v^{(0)}}$ and we obtain a differential series $G(z,\ldots,D_l^nz)$ with Weierstrass order 1. Moreover \begin{center}
$\textrm{Supp}\ G=\textrm{Supp}\ \tilde{F}_{k,l}-v^{(0)}\subset (\mathcal{R}+\mathbb{N}\sigma^{(h)}-v^{(0)})_{\geq\underline{0}}$.
\end{center} 
The associated series is \begin{center}
$z_0=\tilde{y}_0/t^{\sigma^{(h)}}= (y_0-p_{S''_h})/t^{\sigma^{(h)}}\in\mathds{K}_r^\prec$
\end{center} that has initial part $z_{0,l}\in\mathds{K}_{r,l}$. For any initial segment $S'$ of $\textrm{Supp}\ z_{0,l}$, we denote by $q_{S'}$ the corresponding initial part of $z_{0,l}$ and $\tilde{v}_{S'}=v(G(q_{S'}))$.

From the induction hypothesis, there exists a well-ordered subset $\mathcal{R}'$ of $\Gamma_{>0}$ obtained from $\textrm{Supp}\ G$, $\mathcal{T}_l,\ldots,\mathcal{T}_r$ by finitely many elementary transformations such that:
\begin{description}
    \item[(a)] either $\textrm{Supp}\ z_{0,l}\subset\mathcal{R}'$;
\item[(b)]or $z_0$ stabilises on $G$ with initial part $q_{S'_0}$ where $S'_0$ is a proper initial segment of $\textrm{Supp}\ z_{0,l}$ which is included in $\mathcal{R}'$.
\end{description}
 We note that $\mathcal{R}'$ is obtained from $\textrm{Supp}\ G$, $\mathcal{T}_l,\ldots,\mathcal{T}_r$ by a finite number of elementary transformations and \begin{center}
$\textrm{Supp}\ G\subset (\mathcal{R}+\mathbb{N}\sigma^{(h)}-v^{(0)})_{\geq\underline{0}}$.
\end{center} So, in case (ii), $S'_0\subset\mathcal{R}'$ where $\mathcal{R}'$ is obtained from $\mathcal{R}$, $\textrm{Supp}\ F$, $\mathcal{T}_k,\ldots,\mathcal{T}_r$ by a finite number of elementary transformations. Moreover, we observe that any proper initial segment $S'$ of $\textrm{Supp}\ z_{0,l}$ corresponds to a proper initial segment $S'+\sigma^{(h)}$ of $\textrm{Supp}\ y_{1,k}$, and so to a proper initial segment $S''_h\cup(S'+\sigma^{(h)})$ of $\textrm{Supp}\ y_{0,k}$. Thus, for any initial segment $S'$ of $\textrm{Supp}\ z_{0,l}$, \begin{center}
$\tilde{v}_{S'}=v(G(q_{S'}))=v(\tilde{F}_{k,l}(q_{S'}))+v^{(0)}$.
\end{center}
 But \begin{center}
$v(\tilde{F}_{k,l}(q_{S'}))=v(\hat{F}(q_{S'})) =v(F_{S''_h}(p_{S'+\sigma^{(h)}}))$
\end{center} and \begin{center}
$v(F_{S''_h}(p_{S'+\sigma^{(h)}}))=v(f_{S''_h\cup (S'+\sigma^{(h)})})=v_{S''_h\cup(S'+\sigma^{(h)})}$.
\end{center} So $\tilde{v}_{S'}=v_{S''_h\cup(S'+\sigma^{(h)})}+v^{(0)}$.
It means that, in the case (b), $y_0$ stabilises on $F$ with initial segment $p_{S_0}$ where $S_0=S''_h\cup(S'_0+\sigma^{(h)})$. It suffices then to set \begin{center}
$\mathcal{R}_{h,k}=\left\langle\mathcal{R}'\right\rangle+ \mathbb{N}\sigma^{(h)}$.
\end{center}
In case (a), $\textrm{Supp}\ z_{0,l}\subset\mathcal{R}'$. So \begin{center}
$\textrm{Supp}\ (t^{\sigma^{(h)}}z_{0,l})\subset\mathcal{R}'+\sigma^{(h)}$.
\end{center} We set $\tilde{S}'_h=S'_h\cup\textrm{Supp}\ (t^{\sigma^{(h)}}z_{0,l})$, $\tilde{\sigma}^{(h)}$ its successor element in $\textrm{Supp}\ y_{0,k}$ and $\tilde{S}''_h=\tilde{\sigma}^{(h)}\cup\{\tilde{\sigma}^{(h)}\}$. If we denote $\tilde{z}_0=(y_0-p_{\tilde{S}''_h})/t^{\tilde{\sigma}^{(h)}}$, then $\tilde{z}_0\in\mathds{K}_r^\prec$ with initial part $\tilde{z}_{0,l'}\in\mathds{K}_{r,l'}$ for some $l'\in\{k+1,\ldots,l-1\}$. Therefore we repeat the previous process, starting with the additive conjugation $y=y_1+p_{\tilde{S}''_h}$.

There are at most as many steps as \begin{center}
$z_0=p_{S_h\backslash\{\sigma^{(h)}\}}/t^{\sigma^{(h)}}= z_{0,l}+\cdots+z_{0,k+1}$
\end{center} has elements $z_{0,i}\in\mathds{K}_{r,i}$ in its analysis. If we denote by $j$ this number of steps, in the case where there is no stabilisation, we obtain for any $i=0\cdots j-1$, \begin{center}
$\textrm{Supp}\ (t^{\tilde{\sigma}_{(h)}^{(i)}}\tilde{z}^{(i)}_{0,l^{(i)}}) \subset\mathcal{R}^{(i+1)}$
\end{center} where $\mathcal{R}^{(i+1)}$ is a well-ordered subset of $\Gamma_{>0}$ obtained from $\mathcal{R}$, $\textrm{Supp}\ F$, and $\mathcal{T}_k,\ldots,\mathcal{T}_r$ by finitely many elementary transformations. But \begin{center}
$S_h=\bigcup_{i=0}^{j-1}\textrm{Supp}\ (t^{\tilde{\sigma}_{(h)}^{(i)}}\tilde{z}^{(i)}_{0,l^{(i)}})$.
\end{center} So $S_h\subset\bigcup_{k=1}^j\mathcal{R}^{(k)}$.
Then we set \begin{center}
$\mathcal{R}_{h,k}=\left\langle\bigcup_{k=1}^j\mathcal{R}^{(k)}\right\rangle$ 
\end{center}which is an additive sub-semigroup of $\Gamma_{>0}$ obtained from $\mathcal{R}$, $\textrm{Supp}\ F$, $\mathcal{T}_k,\ldots,\mathcal{T}_r$ by finitely many elementary transformations.

\subsection{The case $w>1$.}\label{sect:w>1}

We consider a Differential Series \ref{eq:serie-diff-D_k} with Weierstrass order $w>1$, together with a series $y_0=y_{0,k}+\cdots+y_{0,1}\in\mathds{K}_r^\prec$. To prove the Lemma \ref{lemme:princ}, we proceed in three steps. 

First, since $w>1$, there exists at least one multi-index $I\in\mathbb{N}^{n+1}$ with $|I|=w-1$ such that $F^{(I)}$ has Weierstrass order 1. Thus we apply the induction hypothesis to these differential series, and obtain the desired finiteness property for (at least) some proper initial part $p_{S_{w-1}}$ of $y_{0,k}$. We denote \begin{center}
$\mu^{(w-1)}=(0,\ldots,0,\mu^{(w-1)}_k,\ldots,\mu^{(w-1)}_r)= v(y_0-p_{S_{w-1}})$ 
\end{center}and $S''_{w-1}=S_{w-1}\cup\{\mu^{(w-1)}\}$.
Then we note that $(y_0-p_{S''_{w-1}})/t^{\mu^{(w-1)}}$ has an initial part $z_{0,l}+\cdots+z_{0,k}$ for some $l\in\{r,\ldots,k\}$.

 In the case where $l\in\{r,\ldots,k+1\}$, the second step is devoted to determine the support of $z_{0,l}+\cdots+z_{0,k+1}$. We reduce $F(y,\ldots,D_k^ny)$ to another series $G(z,\ldots,D_l^nz)$ with Weierstrass order $\tilde{w}$ together with a series associated $z_0$ with initial part $z_{0,l}\in\mathds{K}_{r,l}$ such that $(\tilde{w},l)<(w,k)$, by means of three successive transformations.

In the third step, we check the support of $z_{0,k}$ reducing to a differential series $\tilde{G}(\tilde{z},\ldots,D_k^n\tilde{z})$ of Weierstrass order lower than $w$.\\

\noindent\underline{\textbf{First step}}. We need some new notations.
\begin{notation}We denote: 
\begin{itemize}
    \item $\mathcal{A}=\{I\in\mathbb{N}^{n+1}\ \mid\ |I|=w, v(c_I)=\underline{0}\}$;
\item for any $I\in\mathcal{A},\ c_{I,0}=\delta(c_I)$;
\item $\mathcal{A}_{w-1}=\{I\in\mathbb{N}^{n+1}\ \mid\ \exists J\in\mathbb{N}^{n+1},\ |J|=1,\ I+J\in\mathcal{A}\}$.\\ So we have $|I|=w-1$  for all $I\in\mathcal{A}_{w-1}$.
\item For any given $I\in\mathcal{A}_{w-1}$, we set:
\begin{itemize}
\item $\mathcal{A}_I=\{J\in\mathbb{N}^{n+1}\ \mid\ |J|=1,\ I+J\in\mathcal{A}\}$;
\item $\pi_I(X)=\sum_{J\in\mathcal{A}_I}c_{I+J,0}X^{\|J\|}$ (then we have $\deg \pi_I\leq n$);
\item  $\{\rho_1^{(I)},\ldots,\rho_{m_I}^{(I)}\}$ the set of the positive roots of $\pi_I$.
\end{itemize}
\end{itemize}
\end{notation}

Then we apply the induction hypothesis and the Lemma \ref{lemma:dichot_weier1} to the differential series $F^{(I)}(y,\ldots,D_k^ny)$, $I\in\mathcal{A}_{w-1}$ together with the series $y_0$. We obtain two cases:
\begin{itemize}
\item either there exists $I_0\in\mathcal{A}_{w-1}$ such that we are in the first case of Lemma \ref{lemme:princ}: there exists a well-ordered subset $\mathcal{R}_{I_0}$ of $\Gamma_{>0}$ obtained from $\textrm{Supp}\ F^{(I_0)},\ \mathcal{T}_1,\ldots,\mathcal{T}_r$ by a finite number of elementary transformations such that $\textrm{Supp}\ y_{0,k}\subset\mathcal{R}_{I_0}$. It suffices then to observe that $\textrm{Supp}\ F^{(I_0)}\subset\textrm{Supp}\ F$;
\item or for any $I\in\mathcal{A}_{w-1}$, $y_0$ stabilises on $F^{(I)}$ with an initial part $p_{S_I\cup\{\mu^{(I)}\}}$, $S_I$ being an initial segment of $\textrm{Supp}\ y_{0,k}$ with successor element $\mu^{(I)}$. Moreover the value of stabilisation is at most equal to $\mu^{(I)}$.
\end{itemize} 
In this last case, we set $S_{w-1}=\bigcup_{I\in\mathcal{A}_{w-1}}S_0^{(I)}$, $\mu^{(w-1)}$ its successor element and $S''_{w-1}=S_{w-1}\cup\{\mu^{(w-1)}\}$. So there exists $I_0\in\mathcal{A}_{w-1}$ such that $S_{w-1}=S_0^{(I)}\subset\mathcal{R}_{I_0}$.
We denote \begin{center}
$\mathcal{R}_1=\left\langle\mathcal{R}_{I_0}\cup\{\mu^{(w-1)}\}\right\rangle= \left\langle\mathcal{R}_{I_0}\right\rangle+\mathbb{N}\mu^{(w-1)}$.
\end{center} So $S''_{w-1}\subset\mathcal{R}_1$. For any $I\in\mathcal{A}_{w-1}$, we denote $\hat{v}^{(I)}$ the value of stabilisation of $v_S^{(I)}$, and so $\hat{v}^{(I)}=v_{S''_{w-1}}^{(I)}\leq \mu^{(w-1)}$. We set $\hat{v}_{w-1}=\min\{\hat{v}^{(I)},\ I\in\mathcal{A}_{w-1}\}$. Thus $\hat{v}_{w-1}\leq \mu^{(w-1)}$.

We examine now the valuations of the differential series $F^{(I)}(y,\ldots,D_k^ny)$ for $I\notin\mathcal{A}_{w-1}$, $|I|=w-1$. From the  Taylor expansion formula (\ref{taylor}), for any successor segment $\tilde{S}$ of $S_{w-1}$, \begin{center}
$f^{(I)}_{S_{w-1}\cup\tilde{S}}=f^{(I)}_{S_{w-1}}+\sum_{|J|\geq 1}(f^{(I+J)}_{S_{w-1}}/J!)p_{\tilde{S}}^{(J)}$.
\end{center} But for any multi-index $J$, $v(p_{\tilde{S}}^{(J)})=|J|v(p_{\tilde{S}})=|J|\mu^{(w-1)}$, and for any $J$ with $|J|=1$, $v(f^{(I+J)}_{S_{w-1}})>\underline{0}$ since $I\notin\mathcal{A}_{w-1}$. Then there are two possibilities: \begin{description}
    \item[(i)] either $v_{S''_{w-1}}^{(I)}>\mu^{(w-1)}$. So, from the ultrametric triangular inequality, $v_{S''_{w-1}\cup\tilde{S}}^{(I)}>\mu^{(w-1)}$ for any successor segment $\tilde{S}$ of $S''_{w-1}$;
    \item[(ii)] or $v_{S''_{w-1}}^{(I)}\leq\mu^{(w-1)}$ which implies that\begin{center}
$v_{S''_{w-1}\cup\tilde{S}}^{(I)}=v_{S''_{w-1}}^{(I)}\leq\mu^{(w-1)}$.
\end{center}  Then we denote by $\hat{v}^{(I)}$ this stabilised valuation, by $\mathcal{B}_{w-1}$ the set of multi-indexes $I$ for whom that situation arises, and \begin{center}
$\tilde{v}_{w-1}=\min\{\hat{v}^{(I)}\ \mid\ I\in\mathcal{A}_{w-1}\cup\mathcal{B}_{w-1}\}$.
\end{center} Thus $\tilde{v}_{w-1}\leq\mu^{(w-1)}$.\\
\end{description}

\noindent\underline{\textbf{Second step}}.
\textbf{(1)} We perform the \textsl{additive conjugation} $y=\tilde{y}+p_{S''_{w-1}}$. From Proposition \ref{transl:propos}, we obtain a differential series \begin{center}
$F_{S''_{w-1}}(\tilde{y},\ldots,D_k^n\tilde{y})= \sum_{I\in\mathbb{N}^{n+1}}
c'_I\tilde{y}^{(I)}$
\end{center} where $c'_I=f^{(I)}_{S_{w-1}}/I!$. In particular, there exists $I_0\in\mathbb{N}^{n+1}$ with $|I_0|=w-1$ such that \begin{center}
$v(c'_{I_0})=\tilde{v}_{w-1}\leq\mu^{(w-1)}$.
\end{center} From the Remark \ref{rem:weier_transl}, $F_{S''_{w-1}}(\tilde{y},\ldots,D_k^n\tilde{y})$ has also Weierstrass order $w$. Moreover  \begin{center}
$\textrm{Supp}\ F_{S''_{w-1}}\subset\textrm{Supp}\ F+\left\langle S''_{w-1}\right\rangle+\mathcal{T}_k\subset\textrm{Supp}\ F+\mathcal{T}_k+\mathcal{R}_1$.
\end{center} The associated series is $\tilde{y}_0=y_0-p_{S''_{w-1}}$ with initial part $\tilde{y}_{0,k}=y_{0,k}-p_{S''_{w-1}}$.\\

\noindent\textbf{(2)} We perform the \textsl{multiplicative conjugation} $\tilde{y}=t^{\mu^{(w-1)}}z$ in $F_{S''_{w-1}}$. From the Proposition \ref{eclt:propos}, we obtain a differential series \begin{center}
$\hat{F}(z,\ldots,D_k^nz)=\sum_{J\in\mathbb{N}^{n+1}}
\hat{c}_Jz^{(J)}$
\end{center} with \begin{center}
$\textrm{Supp}\ \hat{F}\subset\textrm{Supp}\ F_{S''_{w-1}}+\mathbb{N}\mu^{(w-1)}+\mathcal{T}_k\subset\textrm{Supp}\ F+\mathcal{T}_k+\mathcal{R}_1$ 
\end{center}
(by Proposition \ref{eclt:propos}, since $\mathbb{N}\mu^{(w-1)}+\mathcal{R}_1=\mathcal{R}_1$). Moreover \begin{center}
$v(\hat{c}_J)\geq\min\{v(c'_I)+v((t^{\mu^{(w-1)}})^{(K)})\}$ 
\end{center}where the minimum is taken over multi-indexes $I\in\mathbb{N}^{n+1}$ such that $|I|=|J|=|K|,\ \|I\|=\|K\|+\|J\|$ and $J,K\leq_{antilex}I$. Furthermore \begin{center}
$(t^{\mu^{(w-1)}})^{(K)}\sim(\mu^{(w-1)}_k)^{\|K\|}t^{|I|\mu^{(w-1)}}$. 
\end{center}In particular, by the Remark \ref{rem:weier_eclt}, $v(\hat{c}_{J_0})=\tilde{v}_{w-1}+(w-1)\mu^{(w-1)}$ where $J_0$ is the greatest multi-index for the anti-lexicographical ordering among the $J$'s such that $v(c'_J)=\tilde{v}_{w-1}$. So $v(\hat{c}_{J_0})\leq w\mu^{(w-1)}$. Moreover, for any $J$ with $|J|\geq w$, \begin{center}
$v(\hat{c}_{J})\geq|J|\mu^{(w-1)}\geq w\mu^{(w-1)}$.
\end{center} So there exists at least one multi-index $I$ with $|I|<w$ such that \begin{center}
$v(\hat{c}_I)=\min\{v(\hat{c}_J)\ \mid\ J\in\mathbb{N}^{n+1}\}$.
\end{center} We denote $v_{min}$ this minimal valuation and $I_0$ some multi-index with $|I_0|$ as low as possible such that $v(\hat{c}_{I_0})=v_{min}$. Then we observe that, dividing $\hat{F}(z,\ldots,D_k^nz)$ by $t^{v_{min}}$, we obtain a differential series $\tilde{F}(z,\ldots,D_k^nz)$ with Weierstrass order equal to $|I_0|$, and so lower than $w$.  The associated series is $z_{0}=\tilde{y}_{0}/t^{\mu^{(w-1)}} \in\mathds{K}_{r}^\prec$.

If $k\in\{r-1,\ldots,1\}$, then $z_0$ has \textsl{a priori} an initial part \begin{center}
$\tilde{y}_{0,k}/\mu^{(w-1)} =z_{0,l}+\cdots+z_{0,k}$ 
\end{center} for some $l\in\{r,\ldots,k+1\}$. We set \begin{center}
$S_1=\{\mu\in\textrm{Supp}\ y_{0,k}\ \mid\ \mu_k=\mu_k^{(w-1)}, \mu\geq\mu^{(w-1)}\}$
\end{center} (so $S_1$ is a successor segment of $S_{w-1}$). We remark that \begin{center}
$S_1=(\bigcup_{i\in\{k+1,\ldots,l\}}\textrm{Supp}\ z_{0,i})+\mu^{(w-1)}$.
\end{center} Now we determine $\textrm{Supp}\ z_{0,l}$.\\

\noindent\textbf{(3)} We perform the \textsl{change of derivation} \begin{center}
$\hat{F}(z,\ldots,D_k^nz)=\tilde{F}_{k,l}(z,\ldots,D_l^nz)$
\end{center} as in Proposition \ref{chgt_cond1}. Then we show as before (proof of Lemma \ref{lemme:princ} in case $w=1$) that there are two cases. Either $z_0$ stabilises on $\tilde{F}_{k,l}$ with initial part 0, which means that $y_0$ stabilises on $F$ with initial part $p_{S''_{w-1}}$. 

Or dividing $\tilde{F}_{k,l}$ by $t^{v^{(0)}}$ where $v^{(0)}$ is the minimum of the valuation of the coefficients of $\tilde{F}_{k,l}$, we obtain a differential series $G(z,\ldots,D_l^nz)$ with Weierstrass order $\tilde{w}\geq 1$. From Remarks \ref{rem:weier_eclt} and \ref{rem:chge_deriv}, $\tilde{w}\leq w$. Moreover from Proposition \ref{supp_chang-deriv0}, \begin{center}
$\textrm{Supp}\ \tilde{F}_{k,l}\subset\textrm{Supp}\ \hat{F}+\mathcal{T}_k\subset\textrm{Supp}\ F+\mathcal{T}_k+\mathcal{R}_1$.
\end{center} Thus,\begin{center} $\textrm{Supp}\ G\subset(\textrm{Supp}\ F+\mathcal{T}_k+\mathcal{R}_1)_{\geq v^{(0)}}-v^{(0)}$.
\end{center} We apply the induction hypothesis to $G(z,\ldots,D_l^nz)$ together with $z_0$ that has initial part $z_{0,l}$. There exists a well-ordered subset $\tilde{\mathcal{R}}_2$ of $\Gamma_{>0}$ obtained from $\textrm{Supp}\ G,\ \mathcal{T}_l,\ldots,\mathcal{T}_r$ by a finite number of elementary transformations such that:\begin{itemize}
\item either the exponents of $z_{0,l}$ belong to $\tilde{\mathcal{R}}_2$;
\item or the series $z_0$ stabilises on $G$ with initial part $q_0$ which is also a proper initial part of $z_{0,l}$, and the support of $q_0$ is included in $\tilde{\mathcal{R}}_l$. Then it means that $y_0$ stabilises on $F$ with initial part $p_0=p_{S''_{w-1}}+t^{\mu^{(w-1)}}q_0$ and that $\textrm{Supp}\ p_0\subset\mathcal{R}_1+(\tilde{\mathcal{R}}_2+\mu^{(w-1)})$, which is obtained from $\textrm{Supp}\ F,\ \mathcal{T}_k,\ldots,\mathcal{T}_r$ by a finite number of elementary transformations as desired.
\end{itemize}
If the preceding first case holds, we perform the additive conjugation $z=\tilde{z}+t^{\mu^{(w-1)}}z_{0,l}$ in the differential series $\hat{F}(z,\ldots,D_k^nz)$ defined above. We obtain a differential series $\hat{\hat{F}}(z,\ldots,D_k^nz)$ with coefficients that have also as minimal valuation $v_{min}$ for some multi-index $I_0$ (the least multi-index for anti-lexicographical ordering among the coefficients having valuation $v_{min}$: see Remark \ref{rem:weier_transl}). Moreover \begin{center}
$\textrm{Supp}\ \hat{\hat{F}} \subset\textrm{Supp}\ \hat{F}+\left\langle\tilde{\mathcal{R}}_2+\mu^{(w-1)}\right\rangle+\mathcal{T}_k$ 
\end{center}(see Proposition \ref{transl:propos}) which is obtained from $\textrm{Supp}\ F,\ \mathcal{T}_k,\ldots,\mathcal{T}_r$ by a finite number of elementary transformations. The associated series is $\hat{z}_0=z_0-z_{0,l}$, with initial part $z_{0,l_1}$ for some $l_1\in\{l,\ldots,k\}$. If $l_1\neq k$, we resume the preceding arguments.

Thus, in the case where there is no stabilisation, we prove gradually that for any $i\in\{l,\ldots,k+1\}$, $\textrm{Supp}\ z_{0,i}\subset\mathcal{R}_2$ with $\mathcal{R}_2$ that is a well-ordered subset of $\Gamma_{>0}$ obtained from $\textrm{Supp}\ F$, $\mathcal{T}_k,\ldots,\mathcal{T}_r$ by finitely many elementary transformations. Thus \begin{center}
$S_1=\textrm{Supp}\ (z_{0,l}+\cdots+z_{0,k+1})+\mu^{(w-1)}$
\end{center} is included in $\mathcal{R}_2+\mu^{(w-1)}$.\\

\noindent\underline{\textbf{Third step}}. It remains to examine the support of $z_{0,k}$ so as to obtain the desired property for the set \begin{center}
$\{\mu\in\textrm{Supp}\ y_{0,k}\ \mid\ \mu_k > \mu_k^{(w-1)}\}$.
\end{center} We return to the differential series $\tilde{F}=\hat{F}/t^{v_{\min}}$ defined in the second step, which has Weierstrass order $|I_0|<w$. Then we perform the additive conjugation \begin{center}
$z=\tilde{z}+(z_{0,l}+\cdots+z_{0,k+1})$.
\end{center} By the Proposition \ref{transl:propos} and the Remark \ref{rem:weier_transl}, we obtain a new differential series $\tilde{G}(\tilde{z},\ldots,D_k^n\tilde{z})$ with the same Weierstrass order $|I_0|<w$, together with a series $\tilde{z}_0=z_0-(z_{0,l}+\cdots+z_{0,k+1})$ that has initial part $z_{0,k}\in\mathds{K}_{r,k}$. Moreover, \begin{center}
$\textrm{Supp}\ \tilde{G}\subset\textrm{Supp}\ \tilde{F}+\left\langle S_1-\mu^{(w-1)}\right\rangle+\mathcal{T}_k\subset(\textrm{Supp}\ F+\mathcal{T}_k+\mathcal{R}_1)_{\geq v_{\min}}-v_{\min} +\left\langle\mathcal{R}_2\right\rangle +\mathcal{T}_k$.
\end{center} Then we apply the induction hypothesis to $\tilde{G}(\tilde{z},\ldots,D_k^n\tilde{z})$ together with $\tilde{z}_0$. There exists a well-ordered subset $\mathcal{R}_3$ of $\Gamma_{>0}$ obtained from $\textrm{Supp}\ \tilde{G},\ \mathcal{T}_k,\ldots,\mathcal{T}_r$ by finitely many elementary transformations such that:\begin{itemize}
\item either the exponents of $z_{0,k}$ belong to $\mathcal{R}_3$. So if we set \begin{center}
$\mathcal{R}=\mathcal{R}_1+[(\mathcal{R}_2+\mathcal{R}_3)+\mu^{(w-1)}]$, 
\end{center}then $\textrm{Supp}\ y_{0,k}\subset\mathcal{R}$ which is obtained from $\textrm{Supp}\ F$, $\mathcal{T}_k,\ldots,\mathcal{T}_r$ by a finite number of elementary transformations as desired;
\item or the series $\tilde{z}_0$ stabilises on $\tilde{G}$ with initial part $q_0$ which is also a proper initial part of $z_{0,k}$, and the support of $q_0$ is included in $\mathcal{R}_3$. In this case, it means that $y_0$ stabilises on $F$ with initial part $p_0=p_{S''_{w-1}}+t^{\mu^{(w-1)}}q_0$ which is initial part of $y_{0,k}$.  Moreover, $\textrm{Supp}\ p_0\subset\mathcal{R}$. This concludes the proof of the Lemma \ref{lemme:princ}, and therefore the one of the Theorem \ref{theo_princ}.
\end{itemize}



%
\end{document}